\documentclass{amsart}
\usepackage{amssymb,amsmath}

\usepackage[dvips]{graphicx}

\address{Laboratoire
de Math\'ematiques Pures et Appliqu\'ees
Joseph Liouville, Centre Universitaire de la Mi-Voix, 
Maison de la Recherche Blaise Pascal, 50 rue F.Buisson B.P. 699,
62228 Calais Cedex, France}
\email{Cedric.Lecouvey@lmpa.univ-littoral.fr}

\address{Department of Mathematics and Statistics, State University of New York at Albany, Albany, NY 12222}
\email{lenart@albany.edu}

\thanks{The second author was partially supported by National Science Foundation 
grants DMS-0403029 and DMS-0701044}


\theoremstyle{plain}
\newtheorem{theorem}{Theorem}[subsection]

\newtheorem{conjecture}[theorem]{Conjecture}
\newtheorem{corollary}[theorem]{Corollary}

\newtheorem{lemma}[theorem]{Lemma}

\newtheorem{proposition}[theorem]{Proposition}

\theoremstyle{definition}
\newtheorem{definition}[theorem]{Definition}
\newtheorem{example}[theorem]{Example}

\theoremstyle{remark}
\newtheorem{remarks}[theorem]{Remarks}
\newtheorem{remark}[theorem]{Remark}

\begin{document}
\bibliographystyle{plain}

\title[$q$-analogs of weight multiplicities for Lie superalgebras]{On $q$-analogs of weight multiplicities for the Lie superalgebras $\mathfrak{gl%
}(n,m)$ and $\mathfrak{spo(}2n,M)$}
\date{}
\author{C\'{e}dric Lecouvey and Cristian Lenart}

\begin{abstract}
The paper is devoted to the generalization of Lusztig's $q$-analog of
weight multiplicities to the Lie superalgebras $\mathfrak{gl}(n,m)$ and $\mathfrak{%
spo(}2n,M).$ We define such $q$-analogs $K_{\lambda,\mu}(q)$ for the typical modules and for
the irreducible covariant tensor $\mathfrak{gl}(n,m)$-modules of highest weight $\lambda$. For $\mathfrak{gl}%
(n,m),$ the defined polynomials  have nonnegative integer coefficients
if the weight $\mu $ is dominant. For $\mathfrak{spo(}2n,M)$,
we show that the positivity property holds when $\mu $ is dominant and 
sufficiently far from a specific wall of the fundamental chamber. We
also establish that the $q$-analog associated to an irreducible covariant
tensor $\mathfrak{gl}(n,m)$-module of highest weight $\lambda $ and a dominant
weight $\mu $ is the generating series of a simple statistic on the set of
semistandard hook-tableaux of shape $\lambda $ and weight $\mu .$ This
statistic can be regarded as a super analog of the charge statistic
defined by Lascoux and Sch\"{u}tzenberger.
\end{abstract}

\maketitle

\vspace{-5mm}

\section{Introduction}

There has been considerable interest recently in defining and studying $q$-analogs of various coefficients in the representation theory of semisimple Lie algebras and, more generally, of Kac-Moody algebras. An important such $q$-analog is the one defined by Lusztig \cite{lusscq} for the dimension of the space of weight $\mu$ in the irreducible representation of a semisimple Lie algebra with highest weight $\lambda$. This $q$-analog is usually denoted by $K_{\lambda,\mu}(q)$, and is known as a {Kostka-Foulkes polynomial}. Its importance is highlighted by its occurence in various contexts beside the original definition. For instance, it was shown to be an affine Kazhdan-Lusztig polynomial (and therefore has positive integer coefficients), it gives the expansion of a Hall-Littlewood polynomial in the basis of irreducible characters, it encodes the Brylinski-Kostant filtration of weight spaces, and is closely related to the so-called energy function in the representation theory of affine algebras \cite{LeS}. Other $q$-analogs studied recently include the generalization of Lusztig's $q$-analog to symmetrizable Kac-Moody algebras \cite{viskfp}, and various $q$-analogs of branching coefficients for semisimple Lie algebras (see \cite{lecqbc} and the references therein).

In this paper, we define and study a generalization of Lusztig's $q$-analog of weight multiplicities to the most fundamental Lie superalgebras, namely the general linear superalgebras $\mathfrak{gl}(n,m)$ and the orthosymplectic superalgebras $\mathfrak{spo(}2n,M)$. A substantial interest in Lie superalgebras comes from mathematical physicists, due to various physical interpretations. In his fundamental paper, Kac \cite{Kac} classified the simple finite dimensional Lie superalgebras that are not Lie algebras. He also gave a formula  for the characters of the finite dimensional irreducible representations of these superalgebras which are known as typical. Since Kac's paper, the investigation centered on character formulas for the atypical representations, and on extending Kashiwara's theory of crystals \cite{kascbq} to Lie superalgebras. 

Our $q$-analog of weight multiplicities for the typical representations of Lie superalgebras is based on a natural quantization of Kac's character formula. More precisely, it is based on a $q$-partition function which, in a certain sense, is defined in precisely the same way as its counterpart for semisimple Lie algebras, on which the definition of Lusztig's $q$-analog is based. We also define a $q$-analog for the irreducible covariant tensor modules via a quantization of a character formula due to Berele-Regev and Sergeev. The positivity of these $q$-analogs is proved by reducing it to that of Lusztig's $q$-analog. 

It is well-known that Lusztig's $q$-analog for the general linear groups is expressed combinatorially via the Lascoux and Sch\"{u}tzenberger charge statistic on semistandard Young tableaux \cite{lassuc}. We derive here a similar statistic on the set of
semistandard hook-tableaux introduced in \cite{BKK}, and show that it can be used to express the $q$-analog associated to an irreducible covariant
tensor $\mathfrak{gl}(n,m)$-module.

Our $q$-analog $K_{\lambda,\mu}(q)$ for the orthosymplectic algebras does not have the positivity property in general. Nevertheless, we prove that the positivity  holds in an important special case, which is related to a certain stabilization phenomenon occuring when the dominant weight $\mu$ is sufficiently far from a specific wall of the fundamental chamber. It is interesting to note that the stabilized version of Lusztig's $q$-analog for the orthogonal and symplectic algebras can be computed via a combinatorial algorithm; the reason is a connection with the energy function on affine crystals, which is explained in \cite{LeS}. In general, there is no known combinatorial formula for Lusztig's $q$-analogs in types $B-D$; in particular, there is no generalization of the Lascoux-Sch\"utzenberger charge. This suggests that our $q$-analogs for the orthosymplectic algebras will also be hard to compute combinatorially. 

It would be interesting to see whether any of the structures related to Lusztig's $q$-analog (affine Hecke algebras, Hall-Littlewood polynomials, the Brylinski-Kostant filtration, the energy function) could be defined for Lie superalgebras. If so, we expect them to be related in the same way to our $q$-analog.

{\bf Acknowledgements.} We are grateful to Georgia Benkart, Jonathan Brundan, and Ronald King for helpful discussions, as well as to the anonymous referee for the careful reading of the manuscript and the suggestions.

\section{Background}

We recall in this section, some background on classical Lie superalgebras. 
The reader is referred to \cite{fssdls} and \cite{Kac} for a more detailed exposition.

\subsection{The root systems for the Lie superalgebras $\mathfrak{gl}(n,m)$ and $%
\mathfrak{spo}(2n,M)$}

Let $n$ and $m$ be two positive integers. Let $\mathfrak{g}=\mathfrak{g}_{0}\oplus 
\mathfrak{g}_{1}$ be one of the Lie superalgebras $\mathfrak{gl}(n,m),$ $\mathfrak{spo}%
(2n,2m+1)$, or $\mathfrak{spo}(2n,2m)$ over $\mathbb{C}.\;$Let $\mathfrak{h}$ be a
Cartan subalgebra of $\mathfrak{g}_{0}.\;$The root system of $\mathfrak{g}$ is
graded so that the set of positive roots $\Delta ^{+}=\Delta _{0}^{+}\sqcup
\Delta _{1}^{+}$ is the disjoint union of the set of positive even and
positive odd roots.\ Let $\{\delta _{\overline{n}},\ldots,\delta _{\overline{1}%
},\delta _{1},\ldots,\delta _{m}\}$ be a basis of $\mathfrak{h}^{\ast }.\;$%
We denote by $\langle\:\cdot\: ,\:\cdot\: \rangle$ the inner product defined on $\mathfrak{h}^{\ast }
$ by 
\begin{equation*}
\langle\,\delta _{\overline{\imath}},\:\delta _{\overline{\jmath}}\,\rangle:=\delta _{i,j}\,\text{, }%
\;\;\langle\,\delta _{r},\:\delta _{s}\,\rangle:=-\delta _{r,s}\,,\;\;\text{ and }\;\;\langle\,\delta _{%
\overline{\imath}},\:\delta _{r}\,\rangle:=0\,,
\end{equation*}
where $\delta _{u,v}$ is the usual Kronecker symbol.\ 

Among the different choices of a set of simple roots, one is called distinguished; the corresponding even roots are those of $\mathfrak{g}_0$, and the odd one is the lowest weight corresponding to the action of $\mathfrak{g}_1$ on $\mathfrak{g}_0$. Our definitions of $q$-analogs of weight multiplicities can be used for any choice of simple roots. Nevertheless, the results in this paper are based on the distinguished simple roots or on any set in their orbit under the action of the Weyl group of $\mathfrak{g}_0$. For other choices, minor modifications are required in the statement of some of our results. We will now specify the distinguished simple roots, cf. \cite[List of Tables]{fssdls}. For $\mathfrak{gl}(n,m)$ with $n,m\ge 1$, they are
\[\delta_{\overline{n}}-\delta_{\overline{n-1}}\:,\;\;\ldots\:,\;\;\delta_{\overline{2}}-\delta_{\overline{1}}\:,\;\;\delta_{\overline{1}}-\delta_{{1}}\:,\;\;\delta_1-\delta_2\:,\;\;\ldots\:,\;\;\delta_{m-1}-\delta_m\,.\]
For $\mathfrak{spo}(2n,1)$ with $n\ge 1$, the distinguished simple roots are
\[\delta_{\overline{n}}-\delta_{\overline{n-1}}\:,\;\;\ldots\:,\;\;\delta_{\overline{2}}-\delta_{\overline{1}}\:,\;\;\delta_{\overline{1}}\,.\]
For $\mathfrak{spo}(2n,2m+1)$ with $n,m\ge 1$, the distinguished simple roots are
\[\delta_{\overline{n}}-\delta_{\overline{n-1}}\:,\;\;\ldots\:,\;\;\delta_{\overline{2}}-\delta_{\overline{1}}\:,\;\;\delta_{\overline{1}}-\delta_{{1}}\:,\;\;\delta_1-\delta_2\:,\;\;\ldots\:,\;\;\delta_{m-1}-\delta_m\:,\;\;\delta_m\,.\]
For $\mathfrak{spo}(2n,2)$ with $n\ge 1$, the distinguished simple roots are
\[\delta_1-\delta_{\overline{n}}\:,\;\;\delta_{\overline{n}}-\delta_{\overline{n-1}}\:,\;\;\ldots\:,\;\;\delta_{\overline{2}}-\delta_{\overline{1}}\:,\;\;2\delta_{\overline{1}}\,.\]
For $\mathfrak{spo}(2n,2m)$ with $n\ge 1,\,m\ge 2$, the distinguished simple roots are
\[\delta_{\overline{n}}-\delta_{\overline{n-1}}\:,\;\;\ldots\:,\;\;\delta_{\overline{2}}-\delta_{\overline{1}}\:,\;\;\delta_{\overline{1}}-\delta_{{1}}\:,\;\;\delta_1-\delta_2\:,\;\;\ldots\:,\;\;\delta_{m-1}-\delta_m\:,\;\;\delta_{m-1}+\delta_m\,.\]

We give below the sets
of positive roots for $\mathfrak{gl}(n,m)$ and $\mathfrak{spo}(2n,M)$ corresponding to the distinguished simple roots. 

\smallskip

\begin{equation*}
\begin{tabular}{|l|l|l|}
\hline
$\mathfrak{gl}(n,m)\,,\;\;n,m\ge 1$
\begin{tabular}{l}
\  \\ 
\  \\ 
\ 
\end{tabular}
 & $\Delta _{0}^{+}=\left\{ 
\begin{tabular}{l}
$\delta _{\overline{\jmath}}-\delta _{\overline{\imath}}\mid 1\leq i<j\leq n$ \\ 
$\delta _{r}-\delta _{s}\mid 1\leq r<s\leq m$%
\end{tabular}
\right\} $ & $\Delta _{1}^{+}=\left\{ \delta _{\overline{\imath}}-\delta
_{r}\mid 
\begin{tabular}{l}
$1\leq i\leq n$ \\ 
$1\leq r\leq m$%
\end{tabular}
\right\} $ \\ 
\hline
$\mathfrak{spo}(2n,2m+1)\,,$%
\begin{tabular}{l}
\  \\ 
\  \\ 
\  \\ \hspace{-2.5cm} $n\ge 1,\,m\ge 0$
\  \\ 
\ 
\end{tabular}
& $\Delta _{0}^{+}=\left\{ 
\begin{tabular}{l}
$2\delta _{\overline{\imath}}\mid 1\leq i\leq n$ \\ 
$\delta _{\overline{\jmath}}\pm \delta _{\overline{\imath}}\mid 1\leq i<j\leq n$ \\ 
$\delta _{r}\mid 1\leq r\leq m$ \\ 
$\delta _{r}\pm \delta _{s}\mid 1\leq r<s\leq m$%
\end{tabular}
\right\} $ & $\Delta _{1}^{+}=\left\{ \delta _{\overline{\imath}}\pm \delta
_{r},\delta _{\overline{\imath}}\mid 
\begin{tabular}{l}
$1\leq i\leq n$ \\ 
$1\leq r\leq m$%
\end{tabular}
\right\} $ \\ 
\hline
$\mathfrak{spo}(2n,2)\,,\;\;n\ge 1$%
\begin{tabular}{l}
\  \\ 
\  \\ 
\ 
\end{tabular}
& $\Delta _{0}^{+}=\left\{ 
\begin{tabular}{l}
$2\delta _{\overline{\imath}}\mid 1\leq i\leq n$ \\ 
$\delta _{\overline{\jmath}}\pm \delta _{\overline{\imath}}\mid 1\leq i<j\leq n$%
\end{tabular}
\right\} $ & $\Delta _{1}^{+}=\left\{ \delta_1\pm \delta _{\overline{\imath}}\mid 1\leq i\leq n\right\} $ \\ 
\hline
$\mathfrak{spo}(2n,2m)\,,$%
\begin{tabular}{l}
\  \\ 
\  \\ 
\  \\ \hspace{-1.8cm} \vspace{0.4cm} $n\ge 1,\,m\ge 2$
\ 
\end{tabular}
& $\Delta _{0}^{+}=\left\{ 
\begin{tabular}{l}
$2\delta _{\overline{\imath}}\mid 1\leq i\leq n$ \\ 
$\delta _{\overline{\jmath}}\pm \delta _{\overline{\imath}}\mid 1\leq i<j\leq n$ \\ 
$\delta _{r}\pm \delta _{s}\mid 1\leq r<s\leq m$%
\end{tabular}
\right\} $ & $\Delta _{1}^{+}=\left\{ \delta _{\overline{\imath}}\pm \delta
_{r}\mid 
\begin{tabular}{l}
$1\leq i\leq n$ \\ 
$1\leq r\leq m$%
\end{tabular}
\right\} $ \\ \hline
\end{tabular}
\end{equation*}

\smallskip

Here, as well as throughout the paper, we adopt the convention that a set of elements which depend on some index in an empty range is empty. For instance, in the case of $\mathfrak{spo}(2n,1)$ we have 
\[\Delta _{0}^{+}=\{2\delta _{\overline{\imath}},\,\delta _{\overline{\jmath}}\pm \delta _{\overline{\imath}} \mid 1\leq i<j\leq n\}\,,\;\;\;\;\Delta _{1}^{+}=\{ \delta _{\overline{\imath}}\mid 1\leq i\leq n\}\,.\] 

Denote by $\rho _{+}$ and $\rho _{-}$ the half sum of positive even and
positive odd roots, respectively.\ Then set $\rho :=\rho _{+}-\rho _{-}.$

 For any $\beta \in \mathfrak{h}^{\ast }$ such that $\beta =\beta _{%
\overline{n}}^{(0)}\delta _{\overline{n}}+\cdot \cdot \cdot +\beta _{%
\overline{1}}^{(0)}\delta _{\overline{1}}+\beta _{1}^{(1)}\delta
_{1}+\cdot \cdot \cdot +\beta _{m}^{(1)}\delta _{m}$, we will write $%
\beta =(\beta ^{(0)};\beta ^{(1)})$, where $\beta ^{(0)}=(\beta _{\overline{n}%
}^{(0)},\ldots,\beta _{\overline{1}}^{(0)})$ and $\beta ^{(1)}=(\beta
_{1}^{(1)},\ldots,\beta _{m}^{(1)})$.\ For any $\kappa =(\kappa _{1},\ldots,\kappa _{p})\in \mathbb{Z}^{p},$
define $\left| \kappa \right|:=\kappa _{1}+\cdot \cdot \cdot +\kappa _{p}.$ We denote by $W$ the Weyl group of $%
\mathfrak{g}_{0}$, and write $\ell $ for the corresponding length function.\ For
any $w\in W,$ we set $\varepsilon (w):=(-1)^{\ell (w)}$.\ The Weyl group acts
on $\mathfrak{h}^{\ast }$. As usual, for any $w\in W$ and $\beta \in \mathfrak{h}^{\ast }$, we
denote by $w(\beta )$ the action of the element $w$ on $\beta $.
The dot action of $W$ on $\mathfrak{h}^{\ast }$ is defined by $w\circ \beta
:=w(\beta +\rho )-\rho $.

\subsection{Classical root subsystems\label{subsec_subroot}}

The set of positive even roots $\Delta _{0}^{+}$ is the set of positive roots
for the Lie algebra $\mathfrak{g}_{0}$.\ Moreover, $\mathfrak{g}_{0}$ can be
identified with a direct sum $\mathfrak{g}_{n}\oplus \mathfrak{g}_{m}$ of two 
Lie algebras of classical type with ranks $n$ and $m$. This splitting is shown in the table below for the superalgebras we are interested in; in all cases, we have $n,m\ge 1$.


\begin{equation*}
\begin{tabular}{|c|c|c|c|c|}
\hline
$\mathfrak{g}$ & $\mathfrak{gl}(n,m)$ & $\mathfrak{spo}(2n,1)$ & $\mathfrak{spo}(2n,2m+1)$ & 
$\mathfrak{spo}(2n,2m)$ \\ \hline
$\mathfrak{g}_{0}=\mathfrak{g}_{n}\oplus \mathfrak{g}_{m}$ & $\mathfrak{gl}(n)\oplus \mathfrak{gl%
}(m)$ & $\mathfrak{sp}(2n)$ & $\mathfrak{sp}%
(2n)\oplus\mathfrak{so}(2m+1)$ & $ \mathfrak{sp}(2n)\oplus\mathfrak{so}(2m)$ \\ \hline
\end{tabular}
\end{equation*}


Write $\Delta _{n}^{+}$ and $\Delta _{m}^{+}$ for the sets of positive roots
of $\mathfrak{g}_{n}$ and $\mathfrak{g}_{m}$, respectively.\ Then $\Delta
_{0}^{+}=\Delta _{n}^{+}\sqcup \Delta _{m}^{+}.$ Let $W_{n}$ and $W_{m}$ be
the Weyl groups associated to the root systems of $\mathfrak{g}_{n}$ and $\mathfrak{g%
}_{m}$.\ They contain subgroups $S_{n}$ and $S_{m}$ isomorphic to the
symmetric groups of rank $n$ and $m.$ Write $\varepsilon ^{(0)}$ and $%
\varepsilon ^{(1)}$ for the signatures defined on $W_{n}$ and $W_{m}.$ We
have $W=W_{n}\times W_{m}$, and  $%
\varepsilon (w)=\varepsilon (w^{(0)})\,\varepsilon (w^{(1)})$ for any $w=(w^{(0)};w^{(1)})\in W$. The half sums
of positive roots $\rho_+^{(0)}$ and $\rho_+^{(1)}$ of $\Delta _{n}^{+}$ and $%
\Delta _{m}^{+}$ verify $\rho _{+}=(\rho_+^{(0)};\rho_+^{(1)}).$ We define the
dot action of $W_{n}$ on the dual Cartan subalgebra of $\mathfrak{g}_{n}$ by  $u\circ \eta
^{(0)}:=u\left(\eta ^{(0)}+\rho_+^{(0)}\right)-\rho_+^{(0)}$. The dot action of $W_{m}$ on the dual Cartan subalgebra of 
$\mathfrak{g}_{m}$ 
is defined similarly. 

Let $P_{n}$ and $P_{m}$ be the sets of integral
weights of $\mathfrak{g}_{n}$ and $\mathfrak{g}_{m}$, respectively.\ Write $P_{n}^{+}
$ and $P_{m}^{+}$ for the subsets of dominant integral weights of $P_{n}$
and $P_{m}$.  Let $P(n,m)\subset \mathfrak{h}^{\ast }$ be the set of integral weights of $\mathfrak{g}$, which can be regarded as the Cartesian product of $P_{n}$ and $P_{m}.$ Denote
by $P^{+}(n,m)\subset P(n,m)$ the subset of dominant integral weights of $%
\mathfrak{g}$, namely the set of weights $\lambda =(\lambda ^{(0)};\lambda
^{(1)})\in P(n,m)$ such that $\lambda ^{(0)}\in P_{n}^{+}$ and $\lambda
^{(1)}\in P_{m}^{+}$ (see Remark \ref{domweights}). 

For the sake of completeness, we recall the explicit description of the weights $\lambda ^{(0)}=(\lambda _{\overline{n}}^{(0)},\ldots,\lambda _{\overline{1}}^{(0)})$ in $P_n^+$ and $\lambda ^{(1)}=(\lambda _{1}^{(1)},\ldots,\lambda_{m}^{(1)})$ in $P_m^+$. In the case of $\mathfrak{gl}(n,m)$, we have $\lambda _{\overline{k}}^{(0)}\in \mathbb{Z}$ for $k=1,\ldots,n$, and  $\lambda _{\overline{n}}^{(0)}\geq \ldots \geq \lambda _{\overline{1}}^{(0)}$. The condition for $\lambda ^{(1)}$ is completely similar. In the case of $\mathfrak{spo}(2n,M)$, we have $\lambda _{\overline{k}}^{(0)}\in \mathbb{Z}$ for $k=1,\ldots,n$, and  $\lambda _{\overline{n}}^{(0)}\geq \ldots \geq \lambda _{\overline{1}}^{(0)}\geq 0$; thus, $\lambda ^{(0)}$ is a partition with at most $n$ parts. 
Similarly, $\lambda _{k}^{(1)}$ for $k=1,\ldots,m$ are either all in ${\mathbb Z}$ or all in $\frac{1}{2}+\mathbb{Z}$; we also require $\lambda _{1}^{(1)}\geq \ldots \geq \lambda _{m}^{(1)}\geq 0$ if $\mathfrak{g}_{m}=\mathfrak{so}_{2m+1}$, and $\lambda _{1}^{(1)}\geq \ldots \geq \lambda_{m-1}^{(1)}\geq \left| \lambda _{m}^{(1)}\right| \geq 0$ if $\mathfrak{g}_{m}=%
\mathfrak{so}_{2m}.$

The dominant weight $\lambda \in P^{+}(n,m)$ is called {\em typical} if 
\begin{equation*}
\langle\,\lambda +\rho ,\:\alpha \,\rangle\neq 0\text{ for any odd positive root }\alpha \in 
\overline{\Delta }_{1}^{+}\,,
\end{equation*}
where $\overline{\Delta }_{1}^{+}:=\Delta _{1}^{+}$ for $\mathfrak{g=gl}(n,m)$,
and $\overline{\Delta }_{1}^{+}:=\left\{ \delta _{\overline{\imath}}\pm
\delta _{r}\mid 1\leq i\leq n,1\leq r\leq m\right\} $ for $\mathfrak{g=spo}%
(2n,2m+1)$ and $\mathfrak{g=spo}(2n,2m).$ A dominant weight which is not typical
is called {\em atypical}.

 Let $W_{\mathrm{stab}}$ be the largest subgroup in $W$ which
stabilizes $\Delta _{1}^{+}.$

\begin{remarks} (1) For $\mathfrak{g=gl}(n,m)$ and $\mathfrak{spo}(2n,2)$, the
set $\Delta _{1}^{+}$ is stable under the action of the Weyl group $W$. Therefore, we
have $W_{\mathrm{stab}}=W.$ 

(2) For $\mathfrak{g=spo}(2n,2m+1)$ with $m\ge 0$ and $\mathfrak{g=spo}%
(2n,2m)$ with $m\ge 2$, the set of odd positive roots $\Delta _{1}^{+}$ is not
stable under the action of $W$. In this case, we have $W_{\mathrm{stab}}=S_{n}\times
W_{m}\subset W.$
\end{remarks}

\subsection{Typical representations and character formula}

To each $\lambda \in P^{+}(n,m)$ is associated an irreducible $\mathfrak{g}$%
-module of highest weight $\lambda $ that will be denoted by $V(\lambda )$.\
In the sequel, we will only consider dominant weights for which $V(\lambda )$ is finite dimensional. We write $P_{f}^{+}(n,m)$ for the subset of $P^{+}(n,m)$ consisting of such weights.  For $\mathfrak{g=gl}(n,m),$ we have $P_{f}^{+}(n,m)=P^{+}(n,m)$. 

\begin{remark}\label{domweights}
In \cite{Kac}, the condition that $V(\lambda)$ is finite dimensional is incorporated in the definition of a dominant weight. 
\end{remark}

For $\mathfrak{g=spo}(2n,M)$, the dominant weight $\lambda =(\lambda
^{(0)};\lambda ^{(1)})$ must verify an additional condition in order to guarantee finite dimensionality of $V(\lambda)$. In order to explain it, we use the explicit description  of $\lambda ^{(0)}=(\lambda _{\overline{n}}^{(0)},\ldots,\lambda _{\overline{1}}^{(0)})$ and $\lambda ^{(1)}=(\lambda _{1}^{(1)},\ldots,\lambda_{m}^{(1)})$ in Section \ref{subsec_subroot}. Given this description, the irreducible $\mathfrak{spo}(2n,M)$-module $V(\lambda )$ is finite dimensional if and only if 
\begin{equation}
\lambda _{j}^{(1)}=0\text{ for any }j>\lambda _{\overline{1}}^{(0)}.
\label{fini}
\end{equation}
In particular, (\ref{fini}) is verified when $\lambda _{\overline{1}%
}^{(0)}\geq m.$ 

The module $V(\lambda )$ is called typical when $\lambda $ is typical, and
atypical otherwise. For any integral weight $\mu \in P(n,m),$ write $%
V(\lambda )_{\mu }$ for the weight subspace of weight $\mu $ in $V(\lambda )$%
. Then $K_{\lambda ,\mu }=\dim V(\lambda )_{\mu }$ is finite.\ The character
of $V(\lambda )$ is defined by 
\begin{equation*}
\mathrm{char}\:V(\lambda ):=\sum_{\mu \in P(n,m)}K_{\lambda ,\mu }\,e^{\mu }.
\end{equation*}

 Set 
\begin{equation*}
\nabla :=\frac{\prod_{\alpha \in \Delta _{1}^{+}}(e^{\alpha /2}+e^{-\alpha
/2})}{\prod_{\alpha \in \Delta _{0}^{+}}(e^{\alpha /2}-e^{-\alpha /2})}%
=e^{-\rho }\frac{\prod_{\alpha \in \Delta _{1}^{+}}(1+e^{-\alpha })}{%
\prod_{\alpha \in \Delta _{0}^{+}}(1-e^{-\alpha })}.
\end{equation*}
Kac has proved that there exists an analog of the Weyl character formula
for the typical finite dimensional simple modules $V(\lambda )$.

\begin{theorem}
\label{Th_Wey_For}{\rm \cite{Kac}} Consider $\lambda \in P_{f}^{+}(n,m)$ a typical
dominant weight for $\mathfrak{g}$.\ Then 
\begin{equation*}
\mathrm{char}\:V(\lambda )=\nabla \sum_{w\in W}\varepsilon (w)\,e^{w(\lambda
+\rho )}.
\end{equation*}
\end{theorem}

Now define the partition function $\mathcal{P}$ based on the expansion 
\begin{equation}\label{defp}
\frac{\prod_{\alpha \in \Delta _{1}^{+}}(1+e^{\alpha })}{\prod_{\alpha \in
\Delta _{0}^{+}}(1-e^{\alpha })}=\sum_{\beta \in P(n,m)}\mathcal{P}(\beta
)\,e^{\beta }.
\end{equation}
As a corollary of Theorem \ref{Th_Wey_For}, we deduce that the multiplicities $K_{\lambda ,\mu }$ for a typical weight $\lambda$ can be expressed in terms of
the partition function $\mathcal{P}$, just like in the case of the semisimple Lie
algebras.

\begin{corollary}
\label{Cor_expr_K}Consider $\lambda \in P_{f}^{+}(n,m)$ a typical weight and 
$\mu \in P(n,m).$ Then 
\begin{equation*}
K_{\lambda ,\mu }=\sum_{w\in W}\varepsilon (w)\,\mathcal{P}(w\circ \lambda
-\mu ).
\end{equation*}
\end{corollary}

\begin{remarks} (1) When $\lambda $ is atypical, the character formula
of Theorem \ref{Th_Wey_For} does not hold in general (see \cite{Ki}).\ The problem of
determining the characters of the atypical irreducible $\mathfrak{gl}(n,m)$%
-modules has been intensively addressed in the literature. There exist in
this case numerous character formulas of Weyl type appropriate to some
particular atypical dominant weights \cite{barhyd,balfci,jhkcfs,Hu,hkjocf,kawihw,pasgir,sertai}.\ Unfortunately none of
these formulas is known to hold in full generality.\ A general algorithm for
computing the characters of the atypical irreducible $\mathfrak{gl}(n,m)$%
-modules has been first given by Serganova \cite{serklp,Ser}.\ In \cite{Bru},
Brundan has also linked the computation of these characters to the
determination of the canonical bases of certain $U_{q}(\mathfrak{gl}(\infty ))$%
-modules.

(2) For $\mathfrak{g=gl}(n,m),$ we have $w(\rho
_{-})=\rho _{-}$ for any $w\in W$.\ Hence, for any $\beta \in P(n,m),$ we have $%
w\circ \beta =w(\beta +\rho _{+})-\rho _{+}$; in other words, $\rho $ can be
replaced by $\rho _{+}$ in the dot action of $W$. Similarly, when $\mathfrak{g=spo}(2n,M),$
we have  $w\circ \beta =w(\beta +\rho _{+})-\rho _{+}$ for any $\beta \in P(n,m)$ and $w\in W_{\mathrm{stab}%
}.$ Observe that the last
statement does not hold in general for $w\in W\setminus W_{\mathrm{stab}}.$

(3) By restricting from $\mathfrak{g}$ to $\mathfrak{g}_{0}$, we can see that $K_{\lambda ,\mu }=K_{\lambda ,w(\mu) }$ for any $w\in W$. 

\end{remarks}

Another important consequence of Theorem \ref{Th_Wey_For} is a
branching rule for the restriction of the typical module $V(\lambda )$ with $%
\lambda \in P_{f}^{+}(n,m)$ from $\mathfrak{g}$ to $\mathfrak{g}_{0}.$ For any $%
\gamma \in P^{+}(n,m)$, write $V^{\mathfrak{g}_{0}}(\gamma )$ for the irreducible
finite dimensional $\mathfrak{g}_{0}$-module of highest weight $\gamma $, and set 
\begin{equation*}
m_{\lambda ,\gamma }:=[V(\lambda ):V^{\mathfrak{g}_{0}}(\gamma )]\,;
\end{equation*}
namely, $m_{\lambda ,\mu }^{\mathfrak{g}}$ is the number of irreducible
components isomorphic to $V^{\mathfrak{g}_{0}}(\gamma )$ in the restriction $%
V(\lambda )\!\downarrow_{\mathfrak{g}_{0}}^{\mathfrak{g}}$. Set 
\begin{equation}
\prod_{\alpha \in \Delta _{1}^{+}}(1+e^{\alpha })=\sum_{\kappa \in
P(n,m)}c(\kappa )\,e^{\kappa }.  \label{def_c}
\end{equation}

\begin{proposition}
\label{Prop-m}{\rm (}cf. Proposition {\rm 2.11} in {\rm \cite{Kac}}{\rm )} Consider $\lambda \in
P_{f}^{+}(n,m)$ a typical weight and $\gamma \in P^{+}(n,m)$. Then 
\[
m_{\lambda ,\gamma }=\sum_{w\in W}\varepsilon (w)\,c(w\circ \lambda -\gamma )\,.
\]
\end{proposition}

\subsection{Irreducible covariant tensor modules for $\mathfrak{gl}(n,m)$}

In this paragraph, we consider $\mathfrak{g=gl}(n,m).$  Then the set of positive
odd roots can be written $\Delta _{1}^{+}=\{\alpha _{i,r}\mid 1\leq
i\leq n,1\leq r\leq m\}$, where $\alpha _{i,r}:=\delta _{\overline{\imath}%
}-\delta _{r}.$ The superalgebra $\mathfrak{gl}(n,m)$ admits a natural
 irreducible module $V$ with dimension $m+n$ and highest weight $%
\delta _{\overline{n}};$ this module may be atypical. For any positive integer $k,$ the tensor power $%
V^{\otimes k}$ is completely reducible \cite{sertai}. Moreover the irreducible modules
appearing in its decomposition have highest weights of the form $\lambda
=(\lambda ^{(0)};\lambda ^{(1)})$ such that
\begin{itemize}
\item  $\left| \lambda ^{(0)}\right| +\left| \lambda ^{(1)}\right| =k$;
\item  $\lambda ^{(0)}$ is a partition with $n$ parts (possibly equal to $0$);
\item  both $\lambda ^{(1)}$ and its conjugate $\mu:=(\lambda ^{(1)})'$ have at most $m$ parts;   
\item  $\mu _{1}\leq \lambda _{\min }^{(0)}$, where $\mu =(\mu _{1},\ldots,\mu _{m})$, and $\lambda _{\min }^{(0)}$ is the smallest nonzero part of $\lambda ^{(0)}$ if $\lambda ^{(0)}\neq
\emptyset$, whereas $\lambda _{\min }^{(0)}:=0$ otherwise.
\end{itemize}

 These highest weights can be identified with {\em $(n,m)$-hook Young
diagrams} by associating to the pair $\lambda=(\lambda ^{(0)};\lambda ^{(1)})$ a
diagram $Y(\lambda )$ obtained by juxtaposing the Young diagrams of $\lambda
^{(0)}$ and $\mu =(\lambda ^{(1)})^{\prime }$ as illustrated in the example
below. The combinatorics of hook Young diagrams was first studied in \cite{barhyd}.

\begin{example}
\label{exa_Young}The hook Young diagram $Y(\lambda )$ for $\lambda
^{(0)}=(9,7,5)$ and $\lambda ^{(1)}=(4,3,3,2)$ (thus $%
(n,m)=(3,4)$ and $\mu =(4,4,3,1)$) is indicated below.

\includegraphics[width=0.40\textwidth]{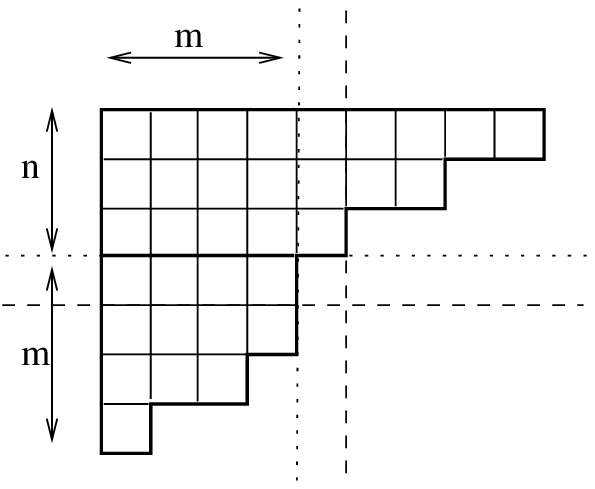}

Observe that there is no box at the intersection of the $(n+1)$-th row and the $%
(m+1)$-th column. By the convention stated earlier, we can write $\lambda =9\delta _{\overline{3}%
}+7\delta _{\overline{2}}+5\delta _{\overline{1}}+4\delta
_{1}+3\delta _{2}+3\delta _{3}+2\delta _{2}$.
\end{example}

 The tensor power $V^{\otimes k}$ contains typical and atypical
irreducible components. Define 
\begin{equation*}
Y_{k}^{+}(n,m):=\left\{ \lambda =(\lambda ^{(0)},\lambda ^{(1)})\mid
V(\lambda )\text{ appears as an irreducible component of a tensor power }%
V^{\otimes k}\right\} .
\end{equation*}
Set $Y^{+}(n,m)=\bigcup _{k\geq 0}Y_{k}^{+}(n,m).$

 Consider $\lambda \in Y^{+}(n,m).$ Set 
\[\nabla _{0}:=e^{-\rho
_{+}}\prod_{\alpha \in \Delta _{0}^{+}}(1-e^{-\alpha })^{-1}\,,\]
 and write $%
\Delta _{1,\lambda }^{+}$ for the subset of $\Delta _{1}^{+}$ containing the
roots $\alpha _{i,r}$ such that $Y(\lambda )$ has a box at the intersection
of its $i$-th row and its $r$-th column. There exists a character formula of
Weyl type for $V(\lambda )$ due to Berele-Regev \cite{barhyd} and Sergeev \cite{sertai}.

\begin{theorem}
\label{prop_char_cova} The character of the covariant tensor module $%
V(\lambda )$ is given by 
\begin{equation*}
\mathrm{char}\:V(\lambda )=\nabla _{0}\sum_{w\in W}\varepsilon (w)\,w\left(
e^{\lambda +\rho _{+}}\prod_{\alpha \in \Delta _{1,\lambda
}^{+}}(1+e^{-\alpha })\right)\,.  
\end{equation*}
\end{theorem}

Set 
\begin{equation}\label{cqlam}\prod_{\alpha \in \Delta _{1,\lambda }^{+}}(1+e^{-\alpha
})=\sum_{\kappa \in P(n,m)}c_{\lambda }(\kappa )\,e^{\kappa }\,.\end{equation}
 We deduce from
the previous proposition the following expression for the multiplicities $%
m_{\lambda ,\gamma }$.

\begin{corollary}\label{mull}
\label{cor_clambda}For any $\lambda \in Y^{+}(n,m)$ and any $\gamma \in
P^{+}(n,m)$, the multiplicity $m_{\lambda ,\gamma }$ of $V^{\mathfrak{g}%
_{0}}(\gamma )$ in $V(\lambda )$ is given by 
\[m_{\lambda ,\gamma }=\sum_{w\in
W}\varepsilon (w)\,c_{\lambda }(\lambda -w\circ \gamma )\,.\]
\end{corollary}

\begin{remarks}\label{remclambda}
(1) When $\Delta _{1,\lambda }^{+}=\Delta _{1}^{+},$
the partition functions $c$ and $c_{\lambda }$ coincide. Moreover, we have $c(w(\kappa
))=c(\kappa )$ for any $\kappa \in P(n,m)$ and any $w\in W.$ Thus the
expression of the multiplicity $m_{\lambda ,\gamma }$  in the 
Corollary \ref{mull} coincides with that in Proposition \ref{Prop-m}.

(2) We cannot derive from the character
formula in Theorem \ref{prop_char_cova} a simple expression for the multiplicity $K_{\lambda
,\mu }$ similar to the one for $\lambda$
typical in Corollary \ref{Cor_expr_K}.
\end{remarks}

\section{Lusztig $q$-analogs and $q$-partition functions}

\subsection{Lusztig $q$-analogs for the Lie algebras $\mathfrak{g}_{n}$ and $\mathfrak{g}_{m}$}

 Denote by $\mathcal{P}_{n,q}$ and $\mathcal{P}_{m,q}$ the $q$%
-partition functions defined by 
\begin{equation}\label{defpqnm}
\frac{1}{\prod_{\alpha \in \Delta _{n}^{+}}(1-qe^{\alpha })}=\sum_{\eta
^{(0)}\in \mathbb{Z}^{n}}\mathcal{P}_{n,q}(\eta ^{(0)})\,e^{\eta ^{(0)}},\;\;\;\;
\frac{1}{\prod_{\alpha \in \Delta _{m}^{+}}(1-qe^{\alpha })}=\sum_{\eta
^{(1)}\in \mathbb{Z}^{m}}\mathcal{P}_{m,q}(\eta ^{(1)})\,e^{\eta ^{(1)}}.
\end{equation}
Consider two weights $\gamma =(\gamma ^{(0)};\gamma ^{(1)})$ and $%
\mu =(\mu ^{(0)};\mu ^{(1)})$ in $P(n,m)$. Set 
\begin{equation}
K_{\gamma ^{(0)},\mu ^{(0)}}^{\mathfrak{g}_{n}}(q):=\sum_{u\in W_{n}}\varepsilon
^{(0)}(u)\,\mathcal{P}_{n,q}(u\circ \gamma ^{(0)}-\mu ^{(0)}),\;\;\;%
K_{\gamma ^{(1)},\mu ^{(1)}}^{\mathfrak{g}_{m}}(q):=\sum_{v\in W_{m}}\varepsilon
^{(1)}(v)\,\mathcal{P}_{m,q}(v\circ \gamma ^{(1)}-\mu ^{(1)}).
\label{def-Lusz}
\end{equation}
These polynomials are Lusztig $q$-analogs for the Lie algebras $\mathfrak{g}%
_{n}$ and $\mathfrak{g}_{m}$.

\begin{theorem}
\label{Th_lusz}{\rm (}Lusztig {\rm \cite{lusscq}}{\rm )} Assume that $\gamma,\mu \in P^{+}(n,m).$ Then
the polynomials $K_{\gamma ^{(0)},\mu ^{(0)}}^{\mathfrak{g}_{n}}(q)$ and $%
K_{\gamma ^{(1)},\mu ^{(1)}}^{\mathfrak{g}_{m}}(q)$ have nonnegative integer
coefficients.
\end{theorem}

We will also need the stabilized form of the
Lusztig $q$-analog corresponding to $\mathfrak{g}_{n}.$ We define
\begin{equation}
K_{\gamma ^{(0)},\mu ^{(0)}}^{\mathfrak{g}_{n},\mathrm{stab}}(q):=\sum_{u\in
S_{n}}\varepsilon ^{(0)}(u)\,\mathcal{P}_{n,q}(u\circ \gamma ^{(0)}-\mu
^{(0)})\,.
\end{equation}
Observe that the sum runs over the symmetric group on $n$ letters.\ For $\mathfrak{%
g=gl}(n,m),$ we have $\mathfrak{g}_{n}\simeq \mathfrak{gl}(n)$ and thus $K_{\gamma
^{(0)},\mu ^{(0)}}^{\mathfrak{g}_{n},\mathrm{stab}}(q)=K_{\gamma ^{(0)},\mu
^{(0)}}^{\mathfrak{g}_{n}}(q).$ For $\mathfrak{g=spo}(2n,M),$ we have $\mathfrak{g}%
_{n}\simeq \mathfrak{sp}_{2n}.\;$Write $\omega $ for the $n$-th fundamental
weight of $\mathfrak{g}_{n}.\;$We can identify $\omega $ with the weight $%
(1,\ldots,1;0,\ldots,0)\in P^{+}(n,m)$. The following lemma has been proved in 
\cite{LeS}.

\begin{lemma}
\label{sta_Cn}For any integer $k\geq \frac{\left| \gamma ^{(0)}\right|
-\left| \mu ^{(0)}\right| }{2}$, we have $K_{\gamma ^{(0)}+k\omega ,\mu
^{(0)}+k\omega }^{\mathfrak{g}_{n}}(q)=K_{\gamma ^{(0)},\mu ^{(0)}}^{\mathfrak{g}%
_{n},\mathrm{stab}}(q).$ In particular, $K_{\gamma ^{(0)},\mu ^{(0)}}^{\mathfrak{g%
}_{n},\mathrm{stab}}(q)\in \mathbb{Z}_{\ge 0}[q]$ when $\gamma ^{(0)},\mu ^{(0)}\in P_{n}^{+}$.
\end{lemma}

\subsection{The Lusztig $q$-analog for $\mathfrak{g}_{0}\label{subsec_Kg0}$}

We define the $q$-partition function $\mathcal{F}_q$ similarly to  (\ref{defpqnm}), by 
\begin{equation}
\frac{1}{\prod_{\alpha \in
\Delta _{0}^{+}}(1-qe^{\alpha })}=\sum_{\eta \in P(n,m)}\mathcal{F}_{q}(\eta
)\,e^{\eta }\,.  \label{Fq}
\end{equation}

\begin{definition}\label{def_Km*n}
For any $\gamma ,\mu \in P(n,m)$, the polynomial $K_{\gamma ,\mu
}^{\mathfrak{g}_{0}}(q)$ is given by
\[
K_{\gamma ,\mu }^{\mathfrak{g}_{0}}(q):=\sum_{w\in W}\varepsilon (w)\,\mathcal{F}%
_{q}(w(\gamma +\rho _{+})-\mu -\rho _{+})\,.  
\]
The
stabilized polynomial $K_{\gamma ,\mu }^{\mathfrak{g}_{0},\mathrm{stab}}(q)$ is given by
\begin{equation*}
K_{\gamma ,\mu }^{\mathfrak{g}_{0},\mathrm{stab}}(q):=\sum_{w\in W_{\mathrm{stab}%
}}\varepsilon (w)\,\mathcal{F}_{q}(w(\gamma +\rho _{+})-\mu -\rho _{+})\,.
\end{equation*}
Since $%
w(\rho _{-})=\rho _{-}$ for any $w\in W_{\mathrm{stab}}$, one can also
replace $\rho _{+}$ by $\rho $ in the latter formula.
\end{definition}

For a dominant weight $\gamma$, the polynomial $K_{\gamma ,\mu }^{\mathfrak{g}_{0}}(q)$ is a $q$-analog for the
dimension $K_{\gamma ,\mu }^{\mathfrak{g}_{0}}$ of the weight space $\mu $ in
the $\mathfrak{g}_{0}$-module of highest weight $\gamma $. From the
considerations in Section \ref{subsec_subroot}, Theorem \ref{Th_lusz}, and Lemma \ref{sta_Cn}, we
deduce easily the following proposition.

\begin{proposition}
\label{prop_U}For $\gamma =(\gamma ^{(0)};\gamma ^{(1)})$ and $\mu =(\mu
^{(0)};\mu ^{(1)})$ in $P(n,m),$ we have, with the above notation, the
factorizations 
\begin{equation*}
K_{\gamma ,\mu }^{\mathfrak{g}_{0}}(q)=K_{\gamma ^{(0)},\mu ^{(0)}}^{\mathfrak{g}%
_{n}}(q)\times K_{\gamma ^{(1)},\mu ^{(1)}}^{\mathfrak{g}_{m}}(q)\,,\;\;\;\;\;\;\;K_{\gamma ,\mu }^{\mathfrak{g}_{0},\mathrm{stab}%
}(q)=K_{\gamma ^{(0)},\mu ^{(0)}}^{\mathfrak{g}_{n},\mathrm{stab}}(q)\times
K_{\gamma ^{(1)},\mu ^{(1)}}^{\mathfrak{g}_{m}}(q)\,.
\end{equation*}
In particular,  $%
K_{\gamma +\omega ,\mu +\omega }^{\mathfrak{g}_{0},\mathrm{stab}}(q)=K_{\gamma
,\mu }^{\mathfrak{g}_{0},\mathrm{stab}}(q).$ Furthermore, $K_{\gamma ,\mu }^{\mathfrak{g}_{0}}(q)$ and $K_{\gamma ,\mu }^{\mathfrak{g}_{0},\mathrm{stab}}(q)$ belong to $%
\mathbb{Z}_{\ge 0}[q]$ when $\gamma$ and $\mu $ are dominant.
\end{proposition}

Consider $\xi \in P(n,m).$ There exists a straightening procedure
for the polynomials $K_{\xi ,\mu }^{\mathfrak{g}_{0}}(q).$

\begin{lemma}
\label{strilaw}Consider $\mu \in P^{+}(n,m)$ and $\xi \in P(n,m).\;$Then 
\begin{equation*}
K_{\xi ,\mu }^{\mathfrak{g}_{0}}(q)=\left\{ 
\begin{tabular}{l}
$\varepsilon (\tau )\,K_{\gamma ,\mu }^{\mathfrak{g%
}_{0}}(q)$ if $\xi =\tau(\gamma+\rho_+)-\rho_+ $ with $\tau \in W$ and $\gamma \in
P^{+}(n,m)$ \\ 
$0$ otherwise.
\end{tabular}
\right. 
\end{equation*}
In particular the coefficients of $K_{\xi ,\mu }^{\mathfrak{g}_{0}}(q)$ are
integers with the same sign.
\end{lemma}

\begin{proof}
The proof follows easily from the equality $K_{\tau(\gamma+\rho_+)-\rho_+  ,\mu }^{%
\mathfrak{g}_{0}}(q)=\varepsilon (\tau )\,K_{\gamma ,\mu }^{\mathfrak{g}_{0}}(q)$ for
any $\tau \in W$.
\end{proof}

\begin{remark}\label{frakk}  Set 
\begin{equation*}
\mathfrak{K}_{\gamma ,\mu }^{\mathfrak{g}_{0}}(q):=\sum_{w\in W}\varepsilon (w)\,%
\mathcal{F}_{q}(w(\gamma +\rho )-\mu -\rho ).
\end{equation*}
For $\mathfrak{g=gl}(n,m)$ and $\mathfrak{spo}(2n,2)$, $W$ stabilizes the set of
positive odd roots $\Delta _{1}^{+}.$ Thus we have $\mathfrak{K}_{\gamma ,\mu }^{%
\mathfrak{g}_{0}}(q)=K_{\gamma ,\mu }^{\mathfrak{g}_{0}}(q)$. This equality does not
hold for the other superalgebras $\mathfrak{spo}(2n,M)$.\ This will cause
some complications in Section \ref{Sec_grad_ops}.
\end{remark}

\subsection{A $q$-partition function associated to $\mathfrak g$}

We define the following natural $q$-analog, denoted by $\mathcal{P}_q$, of the partition function $\mathcal{P}$ in (\ref{defp}). This $q$-partition function will be used below to define our $q$-analogs of weight multiplicities for the Lie superalgebras $\mathfrak{gl}(n,m)$ and $\mathfrak{spo}(2n,M)$.

\begin{definition}\label{defpq} The $q$-partition function $\mathcal{P}_q$ is given by
\begin{equation*}
\frac{\prod_{\alpha \in \Delta _{1}^{+}}(1+qe^{\alpha })}{\prod_{\alpha \in \Delta
_{0}^{+}}(1-qe^{\alpha })}=\sum_{\beta \in P(n,m)}\mathcal{P}_{q}(\beta
)\,e^{\beta }\,. 
\end{equation*}
\end{definition}

 Let us justify why $\mathcal{P}_q$ is a natural $q$-analog, in the spirit of the $q$-partition functions in (\ref{defpqnm}), due to Lusztig. As Kac pointed out in \cite{Kac}, $\mathcal{P}(\beta)$ counts the number of ways to express $\beta$ as $\sum_{\alpha\in\Delta^+}n_\alpha \alpha$ with $n_\alpha$ in $\{0,1\}$ for $\alpha$ in $\Delta_1^+$, and $n_\alpha$ in ${\mathbb Z}_{\ge 0}$ for $\alpha$ in $\Delta_0^+$. Thus, the sequences $\pi=(n_\alpha)_{\alpha\in\Delta^+}$ are generalizations of the {\em Kostant partitions} relevant to (\ref{defpqnm}). We can define the length of the Kostant partition $\pi$, as usual, by $\ell(\pi):=\sum_{\alpha\in\Delta^+} n_\alpha$. Then $\mathcal{P}_{q}(\beta)=\sum_\pi q^{\ell(\pi)}$, precisely as the $q$-partition functions in (\ref{defpqnm}) can be expressed.

We define the $q$-partition function  $c_q$ similarly to  (\ref{def_c}), by 
\begin{equation}
\prod_{\alpha \in \Delta _{1}^{+}}(1+qe^{\alpha })=\sum_{\kappa \in
P(n,m)}c_{q}(\kappa )\,e^{\kappa }\,.  \label{cq}
\end{equation}
Then, for any $\beta \in P(n,m)$, we have%
\begin{equation}
\mathcal{P}_{q}(\beta )=\sum_{\kappa +\eta =\beta }c_{q}(\kappa )\,\mathcal{F}%
_{q}(\eta )\,.  \label{Lem_1}
\end{equation}

\section{Some $q$-analogs of weight multiplicities for $\mathfrak{gl}(n,m)$}

In this section we consider $\mathfrak{g=gl}(n,m)$. For any dominant weight $%
\lambda $ which is either typical or belongs to $Y^{+}(n,m)$, and any weight $\mu
\in P(n,m)$, we introduce a $q$-analog $K_{\lambda ,\mu }(q)$ for the
weight multiplicity $K_{\lambda ,\mu }$. We prove that this $q$-analog has
nonnegative integer coefficients when $\mu \in P^{+}(n,m)$.

\subsection{Typical modules\label{subsec_typi}}

We start by defining $K_{\lambda ,\mu }(q)$ for typical modules similarly to (\ref{def-Lusz}).

\begin{definition}\label{defK(q)}
Consider  a typical dominant weight $\lambda $ and $\mu \in P(n,m).$ The polynomial $K_{\lambda ,\mu }(q)$ is given by  
\begin{equation*}
K_{\lambda ,\mu }(q):=\sum_{w\in W}\varepsilon (w)\,\mathcal{P}_{q}(w\circ
\lambda -\mu )\,,  
\end{equation*}
where $\mathcal{P}_{q}$ is the $q$-partition function in Definition {\rm \ref{defpq}}. 
\end{definition}

According to Corollary \ref{Cor_expr_K}, we have $K_{\lambda ,\mu
}(1)=K_{\lambda ,\mu }$, so $K_{\lambda ,\mu }(q)$ is a $q$-analog for
the multiplicity of $\mu $ in $V(\lambda ).$ We clearly have $K_{\lambda
,\mu }(q)\in \mathbb{Z}[q].$

 Now observe that for any $\alpha =(\alpha ^{(0)};\alpha ^{(1)})\in
\Delta _{1}^{+}$, we have $\left| \alpha ^{(0)}\right| =1$.\ This implies
that if $\kappa =(\kappa ^{(0)};\kappa
^{(1)})\in P(n,m)$ is expressed as a sum of distinct positive odd roots, then the number of summands is equal to $0$ 
or $\left| \kappa ^{(0)}\right| .$ Thus, with the notation in (\ref{def_c})
and (\ref{cq}), we obtain $c_{q}(\kappa )=q^{\left| \kappa ^{(0)}\right|
}c(\kappa )$.

 By the definition (\ref{Fq}) of the partition function $\mathcal{F}%
_{q},$ we can have $%
\mathcal{F}_{q}(\eta )\neq 0$ for some $\eta =(\eta ^{(0)};\eta ^{(1)})\in P_{n,m}$  only if $\left| \eta ^{(0)}\right| =\left|
\eta ^{(1)}\right| =0$.\ Based on the previous observations, (%
\ref{Lem_1}) can be rewritten 
\begin{equation}
\mathcal{P}_{q}(\beta )=q^{\left| \beta ^{(0)}\right| }\sum_{\kappa \in
P(n,m)}c(\kappa )\,\mathcal{F}_{q}(\beta -\kappa )\;\;\;\text{ for any }\beta
=(\beta ^{(0)};\beta ^{(1)})\in P(n,m).  \label{dec-Pq_gl}
\end{equation}
The following theorem can be regarded as an analog of Theorem \ref{Th_lusz}
for the Lie superalgebra $\mathfrak{gl}(n,m).$

\begin{theorem}
\label{TH_KGL}Consider a typical dominant weight $\lambda $  and $\mu $ in $%
P(n,m).$ Then 
\begin{equation}
K_{\lambda ,\mu }(q)=q^{\left| \lambda ^{(0)}\right|-\left|\mu ^{(0)}\right|
}\sum_{\gamma \in P^{+}(n,m)}m_{\lambda ,\gamma }\,K_{\gamma ^{(0)},\mu
^{(0)}}(q)\,K_{\gamma ^{(1)},\mu ^{(1)}}(q)\,.  \label{decA}
\end{equation}
In particular, $K_{\lambda ,\mu }(q)$ belongs to $\mathbb{Z}_{\ge 0}[q]$ when $\mu $
is dominant.
\end{theorem}

\begin{proof}
We derive from (\ref{dec-Pq_gl}) the equality 
\begin{equation*}
K_{\lambda ,\mu }(q)=\sum_{w\in W}\varepsilon (w)\sum_{\kappa \in
P(n,m)}c(\kappa )\,q^{\left| \beta ^{(0)}\right| }\,\mathcal{F}_{q}(w(\lambda
+\rho )-(\mu +\kappa +\rho ))\,,
\end{equation*}
where $\beta =w(\lambda +\rho )-(\mu +\rho )$ in the second sum.\ This
notably implies that $\left| \beta ^{(0)}\right| =\left| \lambda
^{(0)}\right| -\left| \mu ^{(0)}\right| $, since $W=S_{n}\times S_{m}$.\ So
we obtain 
\begin{equation*}
K_{\lambda ,\mu }(q)=q^{\left| \lambda ^{(0)}\right| -\left| \mu
^{(0)}\right| }\sum_{w\in W}\varepsilon (w)\sum_{\kappa \in P(n,m)}c(\kappa )\,%
\mathcal{F}_{q}(w(\lambda +\rho -w^{-1}(\kappa ))-(\mu +\rho ))\,.
\end{equation*}
For any $w\in W$ and any $\kappa \in P(n,m)$, we have $c(\kappa )=c(w(\kappa
))$ because $\Delta _{1}^{+}$ is stable under the action of the Weyl group $%
W.$ Thus, by setting $\xi:=w^{-1}(\kappa )$ in previous expression, we can
write 
\begin{align*}
K_{\lambda ,\mu }(q)&=q^{\left| \lambda ^{(0)}\right| -\left| \mu
^{(0)}\right| }\sum_{w\in W}\varepsilon (w)\sum_{\xi \in P(n,m)}c(\xi )\,%
\mathcal{F}_{q}(w(\lambda +\rho -\xi )-(\mu +\rho ))  \label{equK} \\
&=q^{\left| \lambda ^{(0)}\right| -\left| \mu ^{(0)}\right| }\sum_{\xi \in
P(n,m)}c(\xi )\,K_{\lambda -\xi ,\mu }^{\mathfrak{g}_{0}}(q)\,,  \notag
\end{align*}
where the polynomials $K_{\lambda -\xi ,\mu }^{\mathfrak{g}_{0}}(q)$ are those
in Definition \ref{def_Km*n}.\ Now by Lemma \ref{strilaw}, we have $K_{\lambda -\xi
,\mu }^{\mathfrak{g}_{0}}(q)=0$, or there exists $w\in W$ and $\gamma \in
P^{+}(n,m)$ such that $\gamma =w^{-1}\circ (\lambda -\xi ).$ In the latter case, we have $%
\xi =\lambda +\rho -w(\gamma +\rho )$, and it follows that 
\begin{equation*}
K_{\lambda ,\mu }(q)=q^{\left| \lambda ^{(0)}\right| -\left| \mu
^{(0)}\right| }\sum_{w\in W}\varepsilon (w)\sum_{\gamma \in
P^{+}(n,m)}c(\lambda +\rho -w(\gamma +\rho ))\,K_{\gamma ,\mu }^{\mathfrak{g}%
_{0}}(q)\,.
\end{equation*}
Since $c(\xi )=c(w^{-1}(\xi ))$ for any $w\in W$ and $\xi \in P(n,m)$, we
obtain 
\begin{equation*}
K_{\lambda ,\mu }(q)=q^{\left| \lambda ^{(0)}\right| -\left| \mu
^{(0)}\right| }\sum_{\gamma \in P^{+}(n,m)}\sum_{w\in W}\varepsilon
(w)\,c(w\circ \lambda -\gamma )\,K_{\gamma ,\mu }^{\mathfrak{g}_{0}}(q)\,.
\end{equation*}
Now by Proposition \ref{Prop-m}, we derive 
\begin{equation*}
K_{\lambda ,\mu }(q)=q^{\left| \lambda ^{(0)}\right| -\left| \mu
^{(0)}\right| }\sum_{\gamma \in P^{+}(n,m)}m_{\lambda ,\gamma }\,K_{\gamma
,\mu }^{\mathfrak{g}_{0}}(q).
\end{equation*}
Finally, the desired equality and the positivity for $\mu $ dominant
immediately follow from Proposition \ref{prop_U}.
\end{proof}

\begin{remark}\label{computation}
Recalling the $q$-partition functions $c_q$ and $\mathcal {F}_q$ defined in (\ref{cq}) and (\ref{Fq}), we note that, in practice, it is useful to calculate $K_{\lambda ,\mu }(q)$ by combining   Definition \ref{defK(q)} with (\ref{Lem_1}), as follows: 
\begin{equation}\label{qanalog1}
K_{\lambda ,\mu }(q)=\sum_{w\in W}\sum_{\kappa\in P(n,m)} \varepsilon(w)\,c_q(\kappa) \,{\mathcal F}_q(w\circ\lambda-\mu-\kappa)\,.
\end{equation}
Recalling the notation in (\ref{defpqnm}), also observe that 
\[{\mathcal F}_q(\eta)={\mathcal P}_{n,q}(\eta^{(0)})\times {\mathcal P}_{m,q}(\eta^{(1)})\,. \]
We implemented a Maple procedure based on this approach, which  is part of a package made available at \texttt{math.albany.edu/math/pers/lenart/} and \texttt{lmpa.univ-littoral.fr/\symbol{126}lecouvey/}. The package also contains procedures which compute the $q$-analogs for the orthosymplectic superalgebras defined in Section \ref{Sec_grad_ops}. 
The partition functions  ${\mathcal P}_{n,q}$ and ${\mathcal P}_{m,q}$ are computed by using {\em Gelfand-Tsetlin patterns} of type $A-D$ \cite{baztpmc} to represent Kostant partitions, and by generating these patterns recursively. An efficient implementation is used in order to reduce the computational complexity.  
\end{remark}

\begin{example} Using our Maple procedure, we computed $K_{(3,1,-2;4,2,-8),(0,0,0;0,0,0)}(q)$ for $\mathfrak{gl}(3,3)$:
\[ 2 q^{22}   + 8 q^{21}   + 22 q^{20}   + 40 q^{19}   + 57 q^{18}   + 61 q^{17}   + 52 q^{16}   + 33 q^{15}   + 16 q^{14}  + 5 q^{13}   + q^{12}  \,.\]
\end{example}

Recall that a polynomial is called  unimodal if the sequence of its coefficients has the corresponding property. Based on our experiments, we make the following conjecture.

\begin{conjecture}   For any typical dominant weight $\lambda$ and $\mu \in P^+(n,m)$, the polynomial $K_{\lambda ,\mu }(q)$ is unimodal.
\end{conjecture}

According to Theorem \ref{Th_Wey_For}, we can also define the graded
character of the typical module $V(\lambda )$ by 
\begin{equation}
\mathrm{char}_q\:V(\lambda ):=\nabla (q)\sum_{w\in W}\varepsilon
(w)\,e^{w(\lambda +\rho )}\,,  \label{def_qchar}
\end{equation}
where 
\[\nabla (q):=e^{-\rho }\frac{\prod_{\alpha \in \Delta _{1}^{+}}(1+qe^{-\alpha
})}{\prod_{\alpha \in \Delta _{0}^{+}}(1-qe^{-\alpha })}\,.\]
The coefficients of the expansion of $\mathrm{char}_q\:V(\lambda )$ in the basis
of formal exponentials also yield $q$-analogs of weight multiplicities. It
is easy to verify that these $q$-analogs coincide with the polynomials $%
K_{\lambda ,\mu }(q)$, namely we have 
\begin{equation*}
\mathrm{char}_q\:V(\lambda )=\sum_{\mu \in P(n,m)}K_{\lambda ,\mu }(q)\,e^{\mu
}\,.
\end{equation*}
By Theorem \ref{TH_KGL}, we then obtain the
following expression for $\mathrm{char}_q\:V(\lambda )$:%
\begin{equation}
\mathrm{char}_q\:V(\lambda )=\sum_{\mu \in P(n,m)}q^{\left| \lambda
^{(0)}\right| -\left| \mu ^{(0)}\right| }\sum_{\gamma \in
P^{+}(n,m)}m_{\lambda ,\gamma }\,K_{\gamma ,\mu }^{\mathfrak{g}_{0}}(q)\,e^{\mu }.
\label{dec_gra_char}
\end{equation}

\subsection{Irreducible covariant tensor modules}

For $\nu $ a typical dominant weight, it is possible to define the $q$%
-analogs $K_{\nu ,\mu }(q)$ directly from Corollary \ref{Cor_expr_K} or by
introducing the graded character $\mathrm{char}_q\:V(\nu ).$ \linebreak Now assume $%
\lambda \in Y^{+}(n,m).$ In this case, we have seen (cf. Remark \ref{remclambda} $\mathrm{(2)%
}$) that there is no analog of
Corollary \ref{Cor_expr_K} for the multiplicities $K_{\lambda ,\mu }$.\ Therefore, we define a graded version of the character in Theorem \ref{prop_char_cova}.

\begin{definition} Let
\begin{equation*}
\mathrm{char}_q\:V(\lambda ):=\nabla _{0}(q)\sum_{w\in W}\varepsilon
(w)\,w\left( e^{\lambda +\rho _{+}}\prod_{\alpha \in \Delta _{1,\lambda
}^{+}}(1+qe^{-\alpha })\right)\,,
\end{equation*}
where 
\[\nabla _{0}(q):=e^{-\rho _{+}}\prod_{\alpha \in \Delta
_{0}^{+}}(1-qe^{-\alpha })^{-1}\,.\]
\end{definition}

Based on the above graded character, we define the $q$-analog in the obvious way.

\begin{definition} The polynomials $K_{\lambda ,\mu }(q)$ are given by
\[\mathrm{char}_q\:V(\lambda
)=\sum_{\mu \in P(n,m)}K_{\lambda ,\mu }(q)\,e^{\mu }\,.\]
\end{definition}

Clearly, $K_{\lambda
,\mu }(1)$ is equal to $K_{\lambda ,\mu }$.

\begin{theorem}
We have 
\begin{equation*}
\mathrm{char}_q\:V(\lambda )=\sum_{\mu \in P(n,m)}\sum_{\gamma \in
P^{+}(n,m)}q^{\left| \lambda ^{(0)}\right| -\left| \mu ^{(0)}\right|
}\,m_{\lambda ,\gamma }\,K_{\gamma ,\mu }^{\mathfrak{g}_{0}}(q)\,e^{\mu }\,.
\end{equation*}
In particular, the polynomials $K_{\lambda ,\mu }(q)$ also verify {\rm (\ref{decA}%
)}, and thus have nonnegative integer coefficients when $\mu \in P^{+}(n,m).$
\end{theorem}

\begin{proof}
The arguments of the proof are close to those used in the proof of Theorem \ref{TH_KGL}. By (\ref{cqlam}), we can write 
\[\prod_{\alpha \in \Delta _{1,\lambda
}^{+}}(1+qe^{-\alpha })=\sum_{\kappa \in P(n,m)}q^{\left| \kappa
^{(0)}\right| }\,c_{\lambda }(\kappa )\,e^{-\kappa }\,.\]
 By the definition (\ref{Fq})
of the $q$-partition function $\mathcal{F}_{q},$ we then derive 
\begin{equation*}
\mathrm{char}_q\:V(\lambda )=\sum_{\xi \in P(n,m)}\sum_{\kappa \in
P(n,m)}\sum_{w\in W}\varepsilon (w)\,q^{\left| \kappa ^{(0)}\right|
}\,c_{\lambda }(\kappa )\,\mathcal{F}_{q}(\xi )\,e^{w(\lambda +\rho _{+}-\kappa
)-\xi -\rho _{+}}.
\end{equation*}
By setting $\mu:=w(\lambda +\rho _{+}-\kappa )-\xi -\rho _{+},$ this yields 
\begin{equation*}
\mathrm{char}_q\:V(\lambda )=\sum_{\kappa \in P(n,m)}q^{\left| \kappa
^{(0)}\right| }\,c_{\lambda }(\kappa )\sum_{\mu \in P(n,m)}\sum_{w\in
W}\varepsilon (w)\,\mathcal{F}_{q}(w(\lambda +\rho _{+}-\kappa )-\mu -\rho
_{+})\,e^{\mu }.
\end{equation*}
Thus, we derive 
\begin{equation*}
\mathrm{char}_q\:V(\lambda )=\sum_{\kappa \in P(n,m)}\sum_{\mu \in
P(n,m)}q^{\left| \kappa ^{(0)}\right| }\,c_{\lambda }(\kappa )\,K_{\lambda
-\kappa ,\mu }^{\mathfrak{g}_{0}}(q)\,e^{\mu }.
\end{equation*}
By Lemma \ref{strilaw}, we have $K_{\lambda -\kappa ,\mu }^{\mathfrak{g}_{0}}(q)=0$, or
there exists $w\in W$ and $\gamma \in P^{+}(n,m)$ such that $\gamma
=w^{-1}(\lambda +\rho _{+}-\kappa )-\rho _{+}.$ Then we have $\kappa
=\lambda +\rho _{+}-w(\gamma +\rho _{+})$.\ In particular, $\left| \kappa
^{(0)}\right| =\left| \lambda ^{(0)}\right| -\left| \mu ^{(0)}\right| $ and
it follows that 
\begin{equation*}
\mathrm{char}_q\:V(\lambda )=\sum_{\mu \in P(n,m)}q^{\left| \lambda
^{(0)}\right| -\left| \mu ^{(0)}\right| }\sum_{\gamma \in
P^{+}(n,m)}\sum_{w\in W}\varepsilon (w)\,c_{\lambda }(\lambda +\rho
_{+}-w(\gamma +\rho _{+}))\,K_{\gamma ,\mu }^{\mathfrak{g}_{0}}(q)\,e^{\mu }.
\end{equation*}
Since $\lambda +\rho _{+}-w(\gamma +\rho _{+})=\lambda -w\circ \gamma ,$
this yields the theorem by using Corollary \ref{cor_clambda}.
\end{proof}

\begin{remark} The methods used in this paragraph to define $q$%
-analogs of weight multiplicities corresponding to the irreducible
covariant tensor modules can be extended to any highest weight $\mathfrak{gl}%
(n,m)$-module $V(\lambda )$ whose character is given by a Weyl-type character
formula of the form 
\begin{equation*}
\mathrm{char}\:V(\lambda )=\nabla _{0}\sum_{w\in W}\varepsilon (w)\,w\left(
e^{\lambda +\rho _{+}}\prod_{\alpha \in S_{1,\lambda }^{+}}(1+e^{-\alpha
})\right)\,,
\end{equation*}
where $S_{1,\lambda }^{+}$ is a subset of $\Delta _{1}^{+}$ depending on $%
\lambda$ (cf. \cite{pasgir}). The polynomials $K_{\lambda ,\mu }(q)$ are then defined as the
coefficients of the expansion of the corresponding graded character in the
basis of formal exponentials. They are also expressed as in (\ref{decA}), and therefore 
have nonnegative coefficients when $\mu \in P^{+}(n,m).$ Such a situation
also occurs, for instance, when $V(\lambda) $ is a {\em singly atypical module}, that is,
when there exists a unique root $\alpha _{\lambda }$ in $\Delta _{1}^{+}$
such that $\langle\,\lambda +\rho ,\:\alpha _{\lambda }\,\rangle=0.$ In this case, we have  $S_{1,\lambda }^{+}=\Delta _{1}^{+}\setminus \{\alpha _{\lambda }\}$ (see \cite{jhkcfs,Hu}).
\end{remark}

\subsection{Charge statistic on semistandard hook-tableaux}

We now briefly recall some background on the quantum superalgebra $U_{q}(\mathfrak{gl%
}(n,m))$ and the notion of a {\em crystal basis} (cf. \cite{kascbq}) introduced by Benkart, Kang and
Kashiwara \cite{BKK}.\ The reader is referred to the latter paper for a complete exposition.
One associates to the Lie superalgebra $\mathfrak{gl}(n,m)$ its quantized
enveloping algebra $U_{q}(\mathfrak{gl}(n,m))$, which is a Hopf algebra. This
algebra possesses an irreducible module $V_{q}$ of
dimension $m+n.$ The tensor powers $V_{q}^{\otimes k}$ for $k\in \mathbb{Z}_{\ge 0}$ are
completely reducible. The irreducible modules appearing in their
decompositions into irreducible components are the highest weight $U_{q}(%
\mathfrak{gl}(n,m))$-modules $V_{q}(\lambda )$ with $\lambda \in Y^{+}(n,m).$ To
each $\lambda \in Y^{+}(n,m)$ is associated a $U_{q}(\mathfrak{gl}(n,m))$%
-crystal $B(\lambda ).$ This is an oriented graph with arrows colored with the set of colors
 $\{\overline{n-1},\ldots,\overline{1},0,1,\ldots,m-1\}.$ The vertices of 
$B(\lambda )$ are labelled by semistandard hook-tableaux of shape $\lambda $
on the totally ordered alphabet $\mathcal{A}_{n,m}:=\{\overline{n}<\cdot
\cdot \cdot <\overline{1}<1<\cdot \cdot \cdot <m\}.$ Here, by a semistandard
hook-tableau of shape $\lambda $, we mean a filling of the Young diagram $%
Y(\lambda )$ (see Example \ref{exa_Young}) with letters of $\mathcal{A}_{n,m}$
subject to the following conditions.

\begin{enumerate}
\item  The letters in each row are increasing from left to right.\
Repetition of barred letters is permitted, but repetition of unbarred letters
is not.

\item  The letters in each column are increasing from top to bottom.\
Repetition of unbarred letters is permitted, but repetition of barred letters
is not.
\end{enumerate}

\begin{example}
\label{exa_hooktab}The following tableau is a semistandard hook-tableau of
shape $\lambda =(\lambda ^{(0)},\lambda ^{(1)})$ for $(n,m)=(3,4),$ $\lambda
^{(0)}=(7,6,4)$, and $\lambda ^{(1)}=(4,3,3,2)$.%
\begin{equation*}
\begin{tabular}{|l|llllll}
\hline
$\mathtt{\bar{3}}$ & $\mathtt{\bar{3}}$ & \multicolumn{1}{|l}{$\mathtt{\bar{2%
}}$} & \multicolumn{1}{|l}{$\mathtt{\bar{1}}$} & \multicolumn{1}{|l}{$%
\mathtt{1}$} & \multicolumn{1}{|l}{$\mathtt{2}$} & \multicolumn{1}{|l|}{$%
\mathtt{3}$} \\ \hline
$\mathtt{\bar{2}}$ & $\mathtt{\bar{2}}$ & \multicolumn{1}{|l}{$\mathtt{\bar{1%
}}$} & \multicolumn{1}{|l}{$\mathtt{2}$} & \multicolumn{1}{|l}{$\mathtt{3}$}
& \multicolumn{1}{|l}{$\mathtt{4}$} & \multicolumn{1}{|l}{} \\ \cline{1-6}
$\mathtt{\bar{1}}$ & $\mathtt{\bar{1}}$ & \multicolumn{1}{|l}{$\mathtt{1}$}
& \multicolumn{1}{|l}{$\mathtt{2}$} & \multicolumn{1}{|l}{} &  &  \\ 
\cline{1-4}
$\mathtt{1}$ & $\mathtt{2}$ & \multicolumn{1}{|l}{$\mathtt{3}$} & 
\multicolumn{1}{|l}{$\mathtt{4}$} & \multicolumn{1}{|l}{} &  &  \\ 
\cline{1-4}
$\mathtt{1}$ & $\mathtt{2}$ & \multicolumn{1}{|l}{$\mathtt{3}$} & 
\multicolumn{1}{|l}{$\mathtt{4}$} & \multicolumn{1}{|l}{} &  &  \\ 
\cline{1-1}\cline{1-4}
$\mathtt{2}$ & $\mathtt{3}$ & \multicolumn{1}{|l}{$\mathtt{4}$} & 
\multicolumn{1}{|l}{} &  &  &  \\ \cline{1-1}\cline{1-3}
$\mathtt{2}$ &  &  &  &  &  &  \\ \cline{1-1}
\end{tabular}
.
\end{equation*}
\end{example}

The crystal $B(\lambda )$ has also the structure of a $U_{q}(\mathfrak{g}_{0})$%
-crystal obtained by deleting the arrows colored $0.$ Write $B%
{{}^\circ}%
(\lambda )$ for the crystal obtained in this way.\ Then, for any $\gamma \in
P^{+}(n,m)$, the multiplicity $m_{\lambda ,\gamma }$ is equal to the number
of connected components in $B%
{{}^\circ}%
(\lambda )$ of highest weight $\gamma $.

 Write $B^{\mathfrak{g}_{0}}(\gamma )$ for the abstract $U_{q}(\mathfrak{g}%
_{0})$-crystal of highest weight $\gamma .$ Since $\mathfrak{g}_{0}=\mathfrak{gl}%
(n)\oplus \mathfrak{gl}(m),$ we have $B^{\mathfrak{g}_{0}}(\gamma )=B^{\mathfrak{g}l_{n}}(\gamma
^{(0)})\times B^{\mathfrak{g}l_{m}}(\gamma ^{(1)}),$ i.e., we obtain the direct product
of the $U_{q}(\mathfrak{g}l_{n})$-crystal $B^{\mathfrak{g}l_{n}}(\gamma ^{(0)})$ and 
the $U_{q}(\mathfrak{g}l_{m})$-crystal $B^{\mathfrak{g}l_{m}}(\gamma ^{(1)}).$ In
particular, the corresponding vertices are labelled by the pairs of tableaux $%
(T^{(0)},T^{(1)})$ such that $T^{(0)}$ (resp. $T^{(1)}$) is semistandard on $%
\{\overline{n},\ldots,\overline{1}\}$ (resp.\ on $\{1,\ldots,m\}$) of shape $%
\gamma ^{(0)}$ (resp.\ $\gamma ^{(1)}).$ In \cite{lassuc}, Lascoux and Sch\"{u}%
tzenberger proved that the Lusztig $q$-analog corresponding to the
general linear group can be expressed as the generating series for a
special statistic $\mathrm{ch}$ on semistandard tableaux called {\em charge}
(see \cite{lltpm} for a complete exposition). By the previous arguments and
Proposition \ref{prop_U}, this implies the following proposition.

\begin{proposition}
Consider $\gamma ,\mu \in P^{+}(n,m).$ Then 
\begin{equation*}
K_{\gamma ,\mu }^{\mathfrak{g}_{0}}(q)=\sum_{(T^{(0)},T^{(1)})\in B^{\mathfrak{g}%
_{0}}(\gamma )_{\mu }}q^{\mathrm{ch}(T^{(0)})+\mathrm{ch}(T^{(1)})}
\end{equation*}
where $B^{\mathfrak{g}_{0}}(\gamma )_{\mu }$ is the set of vertices in $B^{\mathfrak{%
g}_{0}}(\gamma )$ with weight $\mu $.
\end{proposition}

 Now consider $T_{\gamma }\in B%
{{}^\circ}%
(\lambda )$ a highest weight vertex of weight $\gamma .\;$Then, the
connected component $B(T_{\gamma })$ of $B%
{{}^\circ}%
(\lambda )$ containing $T_{\gamma }$ is isomorphic to $B^{\mathfrak{g}%
_{0}}(\gamma )$.\ Let $\theta $ be the corresponding isomorphism. Consider $%
T\in B(T_{\gamma })$ and set $\theta (T)=(T^{(0)},T^{(1)}).\;$It is easy to
check that $T^{(0)}$ is the tableau obtained by deleting the unbarred
letters in $T.$ To obtain $T^{(1)}$, we first delete the barred letters in $%
T$.\ This gives a skew tableau that we can conjugate (i.e., reflect in a diagonal), in order to obtain a skew semistandard tableau $T'$.\ Then $T^{(1)}$ is obtained by
rectifying $T'$ via Sch\"{u}tzenberger's {\em jeu de taquin} (e.g., see \cite{fulyt,lasmp}).

It is then natural to define the charge of the semistandard
hook-tableau $T$ by $\mathrm{ch}(T):=\mathrm{ch}(T^{(0)})+\mathrm{ch}%
(T^{(1)})$, where $T^{(0)}$ and $T^{(1)}$ are obtained by the previous
procedure. Since, the polynomials $K_{\lambda ,\mu }(q)$ with $\lambda \in
Y^{+}(n,m)$ verify (\ref{decA}), we obtain the following theorem.

\begin{theorem}
Let $\lambda \in Y^{+}(n,m)$ and $\mu \in P^{+}(n,m)$ be two dominant
weights. Then 
\begin{equation*}
K_{\lambda ,\mu }(q)=q^{\left| \lambda ^{(0)}\right| -\left| \mu
^{(0)}\right| }\sum_{T\in \mathrm{SSHT}(\lambda )_{\mu }}q^{\mathrm{ch}(T)}\,,
\end{equation*}
where $\mathrm{SSHT}(\lambda )_{\mu }$ is the set of semistandard
hook-tableaux of shape $\lambda $ and weight $\mu $.
\end{theorem}

\section{Some $q$-analogs of weight multiplicities for $\mathfrak{spo}(2n,M)%
\label{Sec_grad_ops}$}

In this section, we assume that $\mathfrak{g=spo}(2n,2m+1)$ or $\mathfrak{g=spo}(2n,2m)$. We
introduce $q$-analogs $K_{\lambda ,\mu }(q)$ for the multiplicities $%
K_{\lambda ,\mu }$ corresponding to a typical
weight $\lambda \in P_{f}^{+}(n,m)$. Although  the family of $q$-analogs $K_{\lambda ,\mu }(q)$ for such $\lambda$ and $\mu \in P^{+}(n,m)$  contains polynomials
with negative coefficients, this family possesses a natural subfamily (the stabilized 
$K_{\lambda ,\mu }(q)$) for which the positivity property holds.

\subsection{The polynomials $K_{\protect\lambda ,\protect\mu }(q)$}

We define the $q$-analog $K_{\lambda ,\mu }(q)$ similarly to Definition \ref{defK(q)}. 

\begin{definition}\label{qanalog}
For any weight $\lambda \in P^{+}(n,m)$ and any $\mu \in P(n,m),$ the polynomial $K_{\lambda ,\mu }(q)$ is given by
\begin{equation*}
K_{\lambda ,\mu }(q):=\sum_{w\in W}\varepsilon (w)\,\mathcal{P}_{q}(w\circ
\lambda -\mu )\,,
\end{equation*}
where  $\mathcal{P}_{q}$ is the $q$-partition function  in  Definition {\rm \ref{defpq}}.\
\end{definition}

When $\lambda \in P_{f}^{+}(n,m)$ is typical, the polynomials $%
K_{\lambda ,\mu }(q)$ coincide with the coefficients appearing in the
expansion of the graded character 
\begin{equation*}
\mathrm{char}_q\:V(\lambda )=\nabla (q)\sum_{w\in W}\varepsilon
(w)\,e^{w(\lambda +\rho )}.
\end{equation*}
We have then $K_{\lambda ,\mu }(1)=K_{\lambda ,\mu }$, and therefore $K_{\lambda ,\mu
}(q)$ is a $q$-analog for the dimension of the weight space $\mu $ in the
finite dimensional module $V(\lambda ).$ Observe that the hypothesis $\mu
\in P_{f}^{+}(n,m)$ does not suffice to guarantee that $K_{\lambda ,\mu
}(q)$ belongs to $\mathbb{Z}_{\ge 0}[q]$, as illustrated by the following example.

\begin{example}\label{ex1} Consider $\mathfrak{spo}(2n,2m+1)$ for $n=1$ and $m=2$. We have $\rho=(-2;1,0)+h$, where $h=(\frac{1}{2};\frac{1}{2},\frac{1}{2})$. We verify that $\lambda=(2;1,1)$ is in $P_f^+(n,m)$, and let $\mu=(0;2,1)$ in $P^+(n,m)$. We calculate $K_{\lambda,\mu}(q)$ based on a formula identical to (\ref{qanalog1}), as shown in the table below. The only weights $\kappa$ for which we have nonzero terms in the right-hand side of (\ref{qanalog1}) are $(1;-1,0)$, $(2;-1,-1)$, and $(2;-1,0)$.  
 The columns labeled $u$ and $v$ correspond to those $w=(u;v)$ in $W=W_n\times W_m$  for which ${\mathcal F}_q(w\circ\lambda-\mu-\kappa)\ne 0$. In accordance with standard notation, $s_{\overline{1}}$ is the generator of $W_n$. Also note that, if we view $w$ as a signed permutation of $1,\ldots,n+m$, then $h_w:=w(h)-h$ is a sequence of $-1$'s and $0$'s, where we have a $-1$ in position $i$ precisely when $w$ contains $\overline{\imath}$.  

\begin{equation*}
\begin{tabular}{|c|c|c|c|c|c|c|}
\hline
$\kappa$ & $c_q(\kappa)$ & $\eta:=w\circ \lambda -\mu-\kappa$ & $u$ & $\varepsilon(u)\,{\mathcal P}_{n,q}(\eta^{(0)})$  & $v$ & $\varepsilon(v)\,{\mathcal P}_{m,q}(\eta^{(1)})$ 
 \\ \hline
$(1;-1,0)$ & $q$ & $w(0;2,1)+(1;-2,-1)+h_w$ & $s_{\overline{1}}$ & $-1$ & $Id$ & $1$
 \\ \hline
$(2;-1,-1)$ & $q^2$ & $w(0;2,1)+(0;-2,0)+h_w$ & $Id$ & $1$ & $Id$ & $q$
 \\ \hline
$(2;-1,0)$ & $q^2$ & $w(0;2,1)+(0;-2,-1)+h_w$ & $Id$ & $1$ & $Id$ & $1$
 \\ \hline
\end{tabular}
\end{equation*}

\noindent Hence, we have $K_{(2;1,1),(0;2,1)}(q)=q^3+q^2-q$. Another example is $K_{(5,4,4;3,2,0),(3,2,1;1,1,0)}(q)$ for $\mathfrak{spo}(2n,2m+1)$ with $n=m=3$, which was computed with our Maple procedure: 
\begin{align*}&3 q^{31}+14 q^{30}+52 q^{29}+148 q^{28}+373 q^{27}+817 q^{26}+1640 q^{25}+3000 q^{24}+5132 q^{23}+8174 q^{22}\\&+12283 q^{21}+17338 q^{20}+23138 q^{19}+28977 q^{18}+34022 q^{17}+36993 q^{16}+36953 q^{15}+33259 q^{14}\\&+26478 q^{13}+18045 q^{12}+10121 q^{11}+4332 q^{10}+1211 q^9+97 q^8-65 q^7-17 q^6+q^5\,.\end{align*}
\end{example}

\subsection{Stabilized $q$-analogs}

We now define the subfamily of the polynomials $K_{\lambda,\mu}(q)$ for which the positivity property holds. 

\begin{definition}\label{def_K(q)stab}
For $\lambda \in P^{+}(n,m)$ and $\mu \in P(n,m)$, the polynomial $K_{\lambda ,\mu }^{\mathrm{stab}}(q)$ is given by
\begin{equation*}
K_{\lambda ,\mu }^{\mathrm{stab}}(q):=\sum_{w\in W_{\mathrm{stab}%
}}\varepsilon (w)\,\mathcal{P}_{q}(w\circ \lambda -\mu )\,.  
\end{equation*}
\end{definition}

Note that the sum in the definition of $K_{\lambda ,\mu }^{\mathrm{stab}}(q)$ runs over $W_{\mathrm{stab}}$, and $\lambda $ is not
assumed to be either typical or in $P_{f}^{+}(n,m).$ Moreover, we have $%
K_{\lambda +\omega ,\mu +\omega }^{\mathrm{stab}}(q)=K_{\lambda ,\mu }^{%
\mathrm{stab}}(q)$, where $\omega =(1,\ldots,1;0,\ldots,0)\in P^{+}(n,m).$ This
justifies our terminology.

\begin{lemma}
\label{Lem_stab} Given $\lambda \in P^{+}(n,m)$ and $\mu \in P(n,m)$, there exists $k_{0}\in \mathbb{Z}_{\ge 0}$ such that, for any
nonnegative integer $k\geq k_{0}$, the following conditions are verified:
\begin{enumerate}
\item  $\lambda +k\omega \in P_{f}^{+}(n,m)$;
\item  $\lambda +k\omega $ is typical;
\item  $K_{\lambda +k\omega ,\mu +k\omega }(q)=K_{\lambda
,\mu }^{\mathrm{stab}}(q).$
\end{enumerate}
\end{lemma}

\begin{proof}
We will assume that $\mathfrak{g=spo}(2n,2m+1).$ The case $\mathfrak{g=spo}(2n,2m)$
is similar.

 (1) Write $\lambda +k\omega =\nu =(\nu ^{(0)};\nu ^{(1)}).$ Then $%
\nu _{\overline{1}}^{(0)}=\lambda _{\overline{1}}^{(0)}+k.$ Thus, by (\ref
{fini}), we have $\nu \in P_{f}^{+}(n,m)$ for any $k$ such that $\lambda _{\overline{%
1}}^{(0)}+k\geq m.$

 (2) We have $\rho =(n-m-\frac{1}{2},\ldots,\frac{1}{2}-m;m-%
\frac{1}{2},\ldots,\frac{1}{2})$ and $\overline{\Delta }_{1}^{+}=\left\{ \delta
_{\overline{\imath}}\pm \delta _{r}\mid 1\leq i\leq n,1\leq r\leq m\right\} .
$ Thus, for any $i\in \{1,\ldots,n\}$ and $r\in \{1,\ldots,m\}$, we have 
\begin{align*}&\langle\,\nu
+\rho ,\:\delta _{\overline{\imath}}-\delta _{r}\,\rangle=\lambda _{\overline{\imath}%
}^{(0)}+\lambda _{r}^{(1)}+n+k-i-r+1\;\;\;\mbox{and}\\
&\langle\,\nu +\rho ,\:\delta _{\overline{\imath}%
}+\delta _{r}\,\rangle=\lambda _{\overline{\imath}}^{(0)}-\lambda
_{r}^{(1)}+n-2m+k-i+r\,.\end{align*}
 It is then possible to choose $k$ sufficiently large
so that $\langle\,\nu +\rho ,\:\delta _{\overline{\imath}}-\delta _{r}\,\rangle>0$ and $\langle\,\nu +\rho,\:\delta _{\overline{\imath}}+\delta _{r}\,\rangle>0$ for all 
positive odd roots $\delta _{\overline{\imath}}-\delta _{r}$ and $\delta _{\overline{\imath}}+\delta _{r}$.

 (3) By the definition of the $q$-partition function $\mathcal{F}_{q},$
we have $\mathcal{F}_{q}(\beta )=0$ when $\left| \beta ^{(0)}\right|
<0.$ Consider $w=(w^{(0)},w^{(1)})\in W$ such that $w\notin W_{\mathrm{stab}%
}.$ Since $w^{(0)}\notin S_{n},$ the signed permutation $w^{(0)}$ changes the sign of at least one
coordinate in $\nu^{(0)} +\rho^{(0)} .$ Thus we have 
\[\left| w(\lambda ^{(0)}+k\omega^{(0)}+\rho
^{(0)} )-(\mu ^{(0)}+k\omega^{(0)} +\rho^{(0)})\right| <\left| \lambda
^{(0)}\right| -\left| \mu ^{(0)}\right| -2k\,.\]
 Then for any integer $k\geq 
\frac{\left| \lambda ^{(0)}\right| -\left| \mu ^{(0)}\right| }{2},$ we
have $\mathcal{F}_{q}(w\circ \nu -\mu )=0$, and therefore $K_{\lambda +k\omega
,\mu +k\omega }(q)=K_{\lambda ,\mu }^{\mathrm{stab}}(q).$
\end{proof}

\begin{remarks} (1) Penkov and Serganova \cite{pasgir} gave a character formula for  a simple Lie superalgebra and a corresponding generic weight $\lambda$. Their definition of a generic weight is analogous to the conditions (1)-(3) in Lemma \ref{Lem_stab}.

(2) By Lemma \ref{Lem_stab}, for any $\lambda \in
P^{+}(n,m)$ and any $\mu \in P(n,m),$ the polynomial $K_{\lambda ,\mu }^{\mathrm{stab}}(q)$
can be regarded as a $q$-analog for the multiplicity of the weight $\mu
+k_{0}\omega $ in the finite dimensional module $V(\lambda +k_{0}\omega ).$

(3) Suppose that $\lambda \in P^{+}(n,m)$ and $\mu
\in P(n,m)$ are such that $K_{\lambda +\omega ,\mu +\omega }(q)=K_{\lambda
,\mu }(q).$ Then for any nonnegative integer $k,$ we have $K_{\lambda
+k\omega ,\mu +k\omega }(q)=K_{\lambda ,\mu }(q).$ Thus $K_{\lambda ,\mu
}(q)=K_{\lambda ,\mu }^{\mathrm{stab}}(q)$.
\end{remarks}

\begin{theorem}\label{posklm}
\label{Th_posit-stab}Consider $\lambda ,\mu \in P^{+}(n,m).$ Then $%
K_{\lambda ,\mu }^{\mathrm{stab}}(q)$ belongs to $\mathbb{Z}_{\ge 0}[q].$
\end{theorem}

The proof uses the following easy lemma.

\begin{lemma}\label{lemmult}
\label{Lem_util2}Consider $\lambda ,\gamma \in P^{+}(n,m)$ and let 
\[m_{\lambda ,\gamma }^{\mathrm{stab}}:=\sum_{w\in W_{\mathrm{stab}%
}}\varepsilon (w)\,c(w(\lambda +\rho _{+})-\gamma -\rho _{+})\,.\]
 Then $%
m_{\lambda ,\gamma }^{\mathrm{stab}}\in \mathbb{Z}_{\ge 0}$.
\end{lemma}

\begin{proof}[Proof of  Lemma {\rm \ref{lemmult}}]
We observe that $c(\kappa )=0$ for any $\kappa \in P(n,m)$ such that $\left|
\kappa ^{(0)}\right| <0.$ By using arguments similar to those in the
proof of Lemma \ref{Lem_stab} (3), there exists a nonnegative integer $k$ such that $m_{\lambda
,\gamma }^{\mathrm{stab}}=m_{\lambda +k\omega ,\gamma +k\omega }$.
\end{proof}

\begin{proof}[Proof of Theorem {\rm \ref{posklm}}]
From the description of $\Delta _{1}^{+}$ and (\ref{Lem_1}), we obtain 
\begin{equation*}
\mathcal{P}_{q}(\beta )=\sum_{\kappa \in P(n,m)}q^{\left| \kappa
^{(0)}\right| }\, c(\kappa )\, \mathcal{F}_{q}(\beta -\kappa )\;\;\;\text{ for any }%
\beta =(\beta ^{(0)};\beta ^{(1)})\in P(n,m)\,.
\end{equation*}
Recall also that $w\circ \beta =w(\beta +\rho _{+})-\rho _{+}$ for any $w\in
W_{\mathrm{stab}}.$ By Definition \ref{def_K(q)stab}, this implies that
\begin{equation*}
K_{\lambda ,\mu }^{\mathrm{stab}}(q)=\sum_{\kappa \in P(n,m)}c(\kappa
)\,q^{\left| \kappa ^{(0)}\right| }\sum_{w\in W_{\mathrm{stab}}}\varepsilon (w)\,%
\mathcal{F}_{q}(w(\lambda +\rho _{+})-(\mu +\kappa +\rho _{+}))\,.
\end{equation*}
Hence 
\begin{equation*}
K_{\lambda ,\mu }^{\mathrm{stab}}(q)=\sum_{w\in W_{\mathrm{stab}%
}}\varepsilon (w)\sum_{\kappa \in P(n,m)}q^{\left| \kappa ^{(0)}\right|
}\,c(\kappa )\,\mathcal{F}_{q}(w(\lambda +\rho _{+}-w^{-1}(\kappa ))-(\mu +\rho
_{+}))\,.
\end{equation*}
Set $\xi =w^{-1}(\kappa )$ in the previous sum. Since $w\in W_{\mathrm{stab}}
$, we have $c(\kappa )=c(\xi )$ and $\left| \kappa ^{(0)}\right| =\left| \xi
^{(0)}\right| $, because $\Delta _{1}^{+}$ is stable under the action of $W_{%
\mathrm{stab}}.$ Thus we can write 
\begin{align*}
K_{\lambda ,\mu }^{\mathrm{stab}}(q)&=\sum_{w\in W_{\mathrm{stab}%
}}\varepsilon (w)\sum_{\xi \in P(n,m)}q^{\left| \xi ^{(0)}\right| }\,c(\xi )\,%
\mathcal{F}_{q}(w(\lambda +\rho _{+}-\xi )-(\mu +\rho _{+})) \\
&=\sum_{\xi \in P(n,m)}q^{\left| \xi ^{(0)}\right| }\,c(\xi )\,K_{\lambda -\xi
,\mu }^{\mathfrak{g}_{0},\mathrm{stab}}(q)\,,
\end{align*}
where $K_{\lambda -\xi ,\mu }^{\mathfrak{g}_{0},\mathrm{stab}}(q)
$ are the polynomials in Definition \ref{def_Km*n}. Now, by Lemma \ref{strilaw}, we have $%
K_{\lambda -\xi ,\mu }^{\mathfrak{g}_{0},\mathrm{stab}}(q)=0$, or there exists $%
w\in W_{\mathrm{stab}}$ and $\gamma \in P^{+}(n,m)$ such that $\gamma
=w^{-1}(\lambda -\xi +\rho _{+})-\rho _{+}.$ Hence, we have $\xi =\lambda
+\rho _{+}-w(\gamma +\rho _{+})$ and $\left| \xi ^{(0)}\right| =\left|
\lambda ^{(0)}\right| -\left| \gamma ^{(0)}\right| .$ Thus 
\begin{equation*}
K_{\lambda ,\mu }^{\mathrm{stab}}(q)=\sum_{w\in W_{\mathrm{stab}%
}}\varepsilon (w)\sum_{\gamma \in P^{+}(n,m)}q^{\left| \lambda ^{(0)}\right|
-\left| \gamma ^{(0)}\right| }\,c(\lambda +\rho _{+}-w(\gamma +\rho
_{+}))\,K_{\gamma ,\mu }^{\mathfrak{g}_{0},\mathrm{stab}}(q)\,.
\end{equation*}
Since $c(\xi )=c(w^{-1}(\xi ))$ for any $w\in W_{\mathrm{stab}}$ and $\xi
\in P(n,m)$, we obtain 
\begin{equation*}
K_{\lambda ,\mu }^{\mathrm{stab}}(q)=\sum_{\gamma \in P^{+}(n,m)}q^{\left|
\lambda ^{(0)}\right| -\left| \gamma ^{(0)}\right| }\sum_{w\in W_{\mathrm{%
stab}}}\varepsilon (w)\,c(w(\lambda +\rho _{+})-\gamma -\rho _{+})\,K_{\gamma
,\mu }^{\mathfrak{g}_{0},\mathrm{stab}}(q)\,.
\end{equation*}
Finally, one derives the expression 
\begin{equation*}
K_{\lambda ,\mu }^{\mathrm{stab}}(q)=\sum_{\gamma \in P^{+}(n,m)}q^{\left|
\lambda ^{(0)}\right| -\left| \gamma ^{(0)}\right| }\,m_{\lambda ,\gamma }^{%
\mathrm{stab}}\,K_{\gamma ,\mu }^{\mathfrak{g}_{0},\mathrm{stab}}(q)\,.
\end{equation*}
By Proposition \ref{prop_U} and Lemma \ref{Lem_util2}, this implies $%
K_{\lambda ,\mu }^{\mathrm{stab}}(q)\in \mathbb{Z}_{\ge 0}[q]$ when $\mu \in
P^{+}(n,m).$
\end{proof}

\begin{corollary}
Consider $\lambda ,\mu \in P^{+}(n,m).$ There exists an integer $k_{0}\in 
\mathbb{Z}_{\ge 0}$ such that, for any integer $k\geq k_{0}$, the
following assertions hold:
\begin{enumerate}
\item  $\lambda +k\omega \in P_{f}^{+}(n,m)$;
\item  $\lambda +k\omega $ is typical;
\item  $K_{\lambda +k\omega ,\mu +k\omega }(q)$ has nonnegative integer
coefficients.
\end{enumerate}
\end{corollary}

\begin{remark} We have already observed that for $\mathfrak{g=spo}%
(2n,2)$ we have $W=W_{\mathrm{stab}}.$ This implies that, for any $\lambda
\in P_{f}^{+}(n,m)$ and any $\mu \in P^{+}(n,m),$ the polynomial $K_{\lambda ,\mu
}(q)=K_{\lambda ,\mu }^{\mathrm{stab}}(q)$ has nonnegative integer
coefficients.
\end{remark}

\begin{example} This is a continuation of Example \ref{ex1}. We first considered $\lambda=(2;1,1)$ and  $\mu=(0;2,1)$ for $\mathfrak{spo}(2n,2m+1)$ with $n=1$ and $m=2$. The smallest $k>0$ for which $\lambda+k\omega$ is typical is $3$. So let us consider $\lambda':=\lambda+3\omega=(5;1,1)$ and $\mu':=\mu+3\omega=(3;2,1)$. It turns out that 
\[K^{\rm stab}_{(2;1,1),(0;2,1)}(q)=K_{(5;1,1),(3;2,1)}(q)=q^3+q^2\,.\]
Indeed, the case corresponding to $\kappa=(1;-1,0)$, which is illustrated on the first row of the table in Example \ref{ex1}, and which produced the negative term $-q$, does not appear anymore. To be more precise, in this case we have $w\circ\lambda'-\mu'-\kappa=w(3;2,1)+(-2;-2,-1)+h_w$; but for $w=(s_{\overline{1}};Id)$  this is $(-6;0,0)$, and clearly ${\mathcal F}_q(-6;0,0)=0$. Next, we considered $\lambda=(5,4,4;3,2,0)$ and $\mu=(3,2,1;1,1,0)$ for $\mathfrak{spo}(2n,2m+1)$ with $n=m=3$. In this case, the smallest $k>0$ for which $\lambda+k\omega$ is typical is $5$. Based on the computations with our package, we have 
\[K^{\rm stab}_{(5,4,4;3,2,0),(3,2,1;1,1,0)}(q)=K_{(10,9,9;3,2,0),(8,7,6;1,1,0)}(q)\,,\]
 and this polynomial is
\begin{align*}&3 q^{31}+14 q^{30}+52 q^{29}+148 q^{28}+373 q^{27}+820 q^{26}+1655 q^{25}+3052 q^{24}+5266 q^{23}+8475 q^{22}\\&+12879 q^{21}+18421 q^{20}+24941 q^{19}+31772 q^{18}+38048 q^{17}+42412 q^{16}+43722 q^{15}+41083 q^{14}\\&+34742 q^{13}+25932 q^{12}+16776 q^{11}+9175 q^{10}+4129 q^9+1476 q^8+395 q^7+70 q^6+6 q^5\,.\end{align*}
\end{example}

As in the case of $\mathfrak{gl}(n,m)$, based on our experiments, we make the following conjecture.

\begin{conjecture}  Given $\lambda ,\mu \in P^{+}(n,m)$, the polynomial $K_{\lambda ,\mu }^{\mathrm{stab}}(q)$ is unimodal.
\end{conjecture}



\end{document}